\documentclass[preprint, 11pt]{article}
\usepackage{amsthm,amsmath,amsfonts,amssymb,mathtools,bm,framed,fixltx2e}
\usepackage{bbm}
\usepackage{mathrsfs}
\usepackage{amscd}
\usepackage{enumerate}
\usepackage[numbers,sort&compress]{natbib}
\usepackage[colorlinks,
linkcolor=blue,
anchorcolor=blue,
citecolor=blue]{hyperref}
\usepackage[top=4cm, bottom=3cm, left=3.2cm, right=3.2cm]{geometry}

\theoremstyle{definition}

\theoremstyle{remark}

\theoremstyle{Hypothesis}

\numberwithin{equation}{section} \numberwithin{lem}{section}
\numberwithin{thm}{section} \numberwithin{prop}{section}
\numberwithin{cor}{section} \numberwithin{rem}{section}\numberwithin{hyp}{section}
\allowdisplaybreaks
\title
{Symmetry groups, fundamental solutions and conservation laws for conformable time fractional partial differential system with variable coefficients
}

\author{Xiaoyu {Cheng}
\thanks{Center for Nonlinear Studies, School of Mathematics, Northwest University, Xi'an, 710127, People's Republic of China,
 Email: xiaoyuchengxy@163.com}
 \and Lizhen Wang
\thanks{Center for Nonlinear Studies, School of Mathematics, Northwest University, Xi'an, 710127, People's Republic of China,
Email: wanglizhen@nwu.edu.cn}
}

\begin{document}
\maketitle
\date{}
\textbf{Abstract}: In this paper, the relationships between Lie symmetry groups and fundamental solutions for a class of conformable time fractional partial differential equations (PDEs) with variable coefficients are investigated. Specifically, the group-invariant solutions to the considered equations are constructed applying symmetry group method and the corresponding fundamental solutions for these systems are established with the help of the above obtained group-invariant solutions and inverting Laplace transformation. In addition, the connections between fundamental solutions for two conformable time fractional  systems are given by equivalence transformation. Furthermore, the conservation laws of these fractional systems are provided using new Noether theorem and obtained Lie algebras.


\medskip\noindent
\section{Introduction}
\qquad Fractional calculus implies that the order of the differential and integral operators are fractional numbers and dates back to 1695 when L'Hospital proposed it in the letter to Leibniz. In the last four decades, fractional calculus has attracted many attentions due to its application in various fields such as science, engineering, physics, biology, etc \cite{IP,BM,LD}.

So far, fractional derivatives have been defined in many forms such as Riemann-Liouville, Caputo and Gr\"{u}nwald-Letnikov operators, among which Riemann-Liouville derivative and Caputo derivative are famous and commonly used in modern research. Compared to classical derivative, they lose some basic properties such as chain rule, Leibniz rule and that Riemann-Liouville fractional derivative of a constant is not zero.
In 2014, a new definition was proposed by Khalil \cite{KHYS} called conformable fractional derivative which depends on the definition of the limit. This kind of derivative is well-behaved and Abdeljawad \cite{AT} proved that the conformable derivative obeys some properties such as chain rule, Gromwall's inequality, fractional Laplace transforms, etc. It was observed in Reference \cite{HS} that conformable derivative reflects the direction and strength of the velocity depending on a function in $(t,\epsilon,\alpha)$ and the rate of change of the function for conformable derivative depends on $\alpha$, which makes it $\alpha$-inclusiveness than classical derivative. In recent years, some other properties of conformable derivative have been studied in References \cite{BOK,CW,ZM,KSJ}.

Symmetry group method, first proposed by Sophous Lie, then extended by Olver \cite{OP}, was proved to be a practical method for the analysis of differential equations. In fact, symmetry group of a differential equation is a group which transforms solutions of the equation to other solutions. Therefore, we can directly utilize the property of  symmetry group to construct complex solutions of a system from trivial solutions. More information about symmetry groups can be found in References \cite{OP,LVO,BA0}.  In recent years, symmetry group method has been extended to deal with fractional differential equations \cite{RAS1} and  verified to be one powerful tool to obtain exact solutions of fractional differential equations, such as (1+1)-dimensional Riemann-Liouville fractional equations \cite{LLZ,LHZ,CHW,CWH}, (2+1)-dimensional Riemann-Liouville fractional equations \cite{RSS,RZ,CWL} and so on. For conformable fractional derivative, Chatibi and et al \cite{CKO} proved that Lie symmetry group can be extended to the conformable differential equation and constructed the formulas of the prolongation of the conformable derivative to obtain the exact solutions of the conformable heat equation.

Conservation law plays an important role in the study of some properties of nonlinear PDEs. The correspondence  between Lie symmetry group and conservation laws of PDEs was given by Noether theorem in Reference \cite{N}. Through Noether theorem, one can also construct conservation laws of differential equations. Recently, by the concept of nonlinear self-adjoint equation, Ibragimov \cite{I,IA} has provided a new conservation theorem to study the conservation laws of arbitrary differential equations.

For scalar PDEs, some scholars proved that Lie symmetry group was closely to the fundamental solution of PDEs. Craddock and his collaborators \cite{CMP} first showed symmetry group method can be used to construct fundamental solutions for PDEs of the form
\begin{equation}\label{equation Car}
u_t=x u_{xx}+f(x)u_x,~~x\geq0,
\end{equation}
when the drift function $f(x)$ satisfies Ricatti equations. Their approaches presented that it is always possible to derive classical integral transforms of fundamental solutions of equation \eqref{equation Car} by symmetry. Furthermore, they improved these results in Reference \cite{CML} and obtained the fundamental solutions of a class of equations of the form \begin{equation}\label{equation Ca}
u_t=\sigma x^\gamma u_{xx}+f(x)u_x-\mu x^r u, ~~\sigma>0,
\end{equation}
where $\gamma,$ $\mu$ and $r$ are constants. In Reference \cite{CM}, they further considered
\begin{equation}\label{equation Ca0}
u_t=\sigma x^\gamma u_{xx}+f(x)u_x-g(x) u, ~~ x\geq0, \sigma>0, \gamma\neq2,
\end{equation}
and proved that nontrivial Lie symmetries yield Laplace transform and Fourier transform of the fundamental solutions of equation \eqref{equation Ca0}.

In addition, Kang and Qu \cite{KQ} developed the approach introduced by Craddock and et al \cite{CM,CML,CMP} to study the relationship between Lie symmetries and fundamental solutions of the following system of parabolic equations with variable coefficients
\begin{equation}\label{KangQu}
\left\{
\begin{aligned}
&u_t=u_{xx}+\frac{c}{x}u_x+mx^kv_x,\\
&v_t=v_{xx}+\frac{c}{x}v_x+nx^ku_x,~~ x>0,
\end{aligned}
\right.
\end{equation}
where $c$, $m$, $n$ and $k$ are constants. They set up certain symmetries admitted by system \eqref{KangQu} and provided the corresponding group invariant solutions based on the obtained symmetries. Moreover, the fundamental solutions of system \eqref{KangQu} were derived from its group invariant solutions by means of inverse Laplace transform.

It has been shown that Lie symmetry group was closely related to the fundamental solution of integer PDEs. A natural question arises: can we construct fundamental solutions of fractional linear PDEs from their symmetries? In some cases, Caputo fractional model and conformable fractional model have similar behavior \cite{PR,PE,EAO}, and sometimes conformable fractional model is even more advantageous such as in tumor-immune interactions \cite{BOK}. More importantly, the construction of the solution of equations in the sense of conformable derivative is easier than Caputo derivative. The time derivative term of  system \eqref{KangQu} can be extended to the time fractional derivative of order $\alpha$ and to the best of our knowledge, there is no literature considering the fundamental solution of time fractional system using symmetry group method. So far, compared with Riemann-Liouville derivative and Caputo derivative for Lie  symmetry group, the conformable fractional derivative is more convenient for calculation. In Reference \cite{XYW}, we studied the following conformable time fractional equation
\begin{equation}\label{equation 1}
\mathcal{T}_t^\alpha u=xu_{xx}+f(x)u_x, ~~ x\geq0,
\end{equation}
where $  0<\alpha\leq 1 $, $\mathcal{T}_t^\alpha $ is the conformable fractional differential operator with order $\alpha$.  Moreover, we constructed the fundamental solutions and conservation laws for equation \eqref{equation 1} based on the obtained symmetries.
In this present paper, we consider the following  conformable time fractional system of parabolic equations
\begin{equation}\label{equation 2}
\left\{
\begin{aligned}
&\mathcal{T}_t^\alpha u=u_{xx}+\frac{c}{x}u_x+mx^kv_x,\\
&\mathcal{T}_t^\alpha v=v_{xx}+\frac{c}{x}v_x+nx^ku_x,~~x>0,
\end{aligned}
\right.
\end{equation}
where $  0<\alpha\leq 1 $. By means of the approaches and formulas utilized in References \cite{CM,CML,CMP,CKO,KQ}, we will compute the fundamental solutions for  system \eqref{equation 2} using  obtained symmetry.

In this paper, we intend to construct fundamental solutions and conservation laws of a class of fractional system with conformable time derivative using Lie symmetries admitted by system \eqref{equation 2}. The definitions and properties related to conformable derivative, Laplace transform and Bessel functions are introduced in Section 2. In Section 3, the fundamental solutions of system \eqref{equation 2} are established and the fundamental solutions for two conformable time fractional  systems can be connected by equivalence transformation. In addition, conservation laws of  system \eqref{equation 2} are obtained in Section 4. At the end of this paper, the concluding remarks are presented in Section 5.

\section{Preliminaries}
\qquad In this Section, we recall some definitions and related properties of conformable fractional calculus, Laplace transform and Bessel function.

\textbf{Definition 2.1} \cite{KHYS} Let $f: [0, \infty) \rightarrow \mathbb R$ and $\alpha\in(0,1]$. The conformable derivative of the function $f(t)$ with order $\alpha$ is defined by
\begin{equation*}\label{conformable derivative0}
\mathcal{T}_t^\alpha (f)(t): =\lim_{\epsilon \rightarrow 0}\frac{f(t+\epsilon t^{1-\alpha})-f(t)}{\epsilon},
\end{equation*}
for all $t>0$.

Next, we provide the following lemma of conformable differential operator.

\textbf{Lemma 2.1} \cite{KHYS} Let $0<\alpha\leq1$ and $f(t)$ be $\alpha$-differentiable at a point $t>0$. Then

(a) $\mathcal{T}_t^\alpha(t^p)=pt^{p-\alpha} ~~\mbox{for all} ~~p\in \mathbb{R},$

(b) $\mathcal{T}_t^\alpha(c)=0,~~ c ~\mbox{is a constant},$

(c) $\mbox{in addition, if } ~f(t) ~\mbox{is differentiable, then} ~\mathcal{T}_t^\alpha (f)(t)=t^{1-\alpha}\frac{\mathrm{d}f}{\mathrm{d}t}.$

\textbf{Definition 2.2} \cite{GPH} Let $z\in \mathbb{C}/(-\infty,0].$
The modified Bessel function $I_\nu(z)$ is given by
\begin{equation*}\label{modified Bessel function}
I_\nu(z)=\sum_{n=0}^\infty\frac{(\frac{z}{2})^{2n+\nu}}{n!\Gamma(\nu+n+1)}.
\end{equation*}

\textbf{Definition 2.3} \cite{GPH} Suppose that a function $f(t)$ is defined in $t\in(0,+\infty)$ and its Laplace transform $\tilde{f}(s)$ is defined by
\begin{equation*}\label{LT}
\tilde{f}(s)=\mathcal{L}(f)(s):=\int_0^{+\infty}f(t)\mathrm{e}^{-st}\mathrm{d}t,
\end{equation*}
and $f(t)$ is called the inverse Laplace transform of $\tilde{f}(s)$ and denoted as $f(t)=\mathcal{L}^{-1} (\tilde{f}(s))(t)$.

We present the following lemma related to the Laplace transform.

\textbf{Lemma 2.2} \cite{CMP,GPH} Let $\mathcal{L}$ denote Laplace transformation in $\lambda$, then it holds that

(i) $\mathcal{L}(\mathrm{e}^{at}f(t))(\lambda)=\tilde{f}(\lambda-a),$ where $a$ is arbitrary constant,

(ii) $\mathcal{L}^{-1}\big(\frac{1}{\lambda^\mu}\mathrm{e}^{\frac{k}{\lambda}}\big)
=\big(\frac{y}{k}\big)^{\frac{\mu-1}{2}}I_{\mu-1}(2\sqrt{ky}),$
where $\mu>0$, $k>0$ and $I_{\mu-1}$ is the modified Bessel function with order $\mu-1$.

\section{Lie symmetry group and fundamental solution for system of PDEs \eqref{equation 2}}
\qquad In this Section, aim at finding out the exact expression of the fundamental solution of system \eqref{equation 2}  based on the Lie algebras admitted by the system. First, we introduce the Lie symmetry group and the definition of fundamental solution of conformable fractional parabolic system.
\subsection{Lie point symmetry and fundamental solution to conformable time fractional parabolic system}
\qquad Consider a conformable time fractional system
\begin{equation}\label{parabolic equations}
\left\{
\begin{aligned}
&\mathcal{T}_t^\alpha u=M(x,t,u,v,u^{(1)},v^{(1)},\cdots,u^{(n)},v^{(n)}),\\
&\mathcal{T}_t^\alpha v=N(x,t,u,v,u^{(1)},v^{(1)},\cdots,u^{(n)},v^{(n)}),
\end{aligned}
\right.
\end{equation}
where $M$ and $N$ are assumed to be smooth in their arguments.
Assume that symmetry group $G$ of system \eqref{parabolic equations} is generated by the following vector field
\begin{equation}\label{vector fields v}
V=\xi(x,t,u,v)\partial_x+\tau(x,t,u,v)\partial_t+\eta(x,t,u,v)\partial_u+\phi(x,t,u,v)\partial_v,
\end{equation}
where $\xi$, $\tau$, $\eta$ and $\phi$ are infinitesimals. And formula of the $n$th prolongation of system \eqref{parabolic equations} is presented as follows
\begin{equation}\label{prolongation of v}
\mathrm{Pr}^{(\alpha,n)}V=V+\eta^{\alpha,t}\partial_{\mathcal{T}_t^\alpha u}+\phi^{\alpha,t}\partial_{\mathcal{T}_t^\alpha v}+\eta^x\partial_{u_x}+\phi^x\partial_{v_x}+\eta^{xx}\partial_{u_{xx}}+\phi^{xx}\partial_{v_{xx}}+\cdots,
\end{equation}
where the formulae of $\eta^x$, $\phi^x$, $\eta^{xx}$, $\phi^{xx}$, $\cdots$ are provided and more detailed and rigorous discussions can be found in References \cite{OP,LVO,BA0}. The expressions for $\eta^{\alpha,t}$, $\phi^{\alpha,t}$ are given as follows
\[
\eta^{\alpha,t}=t^{1-\alpha}\eta_t+t^{1-\alpha}\bigg(\eta_u-\tau_t+\frac{1-\alpha}{t}\tau\bigg)u_t+t^{1-\alpha}(\eta_vv_t-\xi_tu_x-\xi_uu_xu_t-\xi_vu_xv_t-\tau_uu_t^2-\tau_vu_tv_t),\]
\[
\phi^{\alpha,t}=t^{1-\alpha}\phi_t+t^{1-\alpha}\bigg(\phi_v-\tau_t+\frac{1-\alpha}{t}\tau\bigg)v_t+t^{1-\alpha}(\phi_uu_t-\xi_tv_x-\xi_uv_xu_t-\xi_vv_xv_t-\tau_uu_tv_t-\tau_vv_t^2).\]
The vector field $V$ satisfies the following invariant condition
\begin{equation}\label{yantuotiaojian}
\left\{
\begin{aligned}
&\mathrm{Pr}^{(\alpha,n)}V(\mathcal{T}_t^\alpha u-M(x,u^{(n)},v^{(n)})|_{\eqref{parabolic equations}}=0,\\
&\mathrm{Pr}^{(\alpha,n)}V(\mathcal{T}_t^\alpha v-N(x,u^{(n)},v^{(n)})|_{\eqref{parabolic equations}}=0.
\end{aligned}
\right.
\end{equation}
Based on system \eqref{yantuotiaojian}, we can deduce an over-determined system and solve this system to derive the vector field $V$.  To exponentiate a vector field $V$, we solve the following system
\begin{equation}\label{systemodes}
\frac{\mathrm{d}\tilde{x}}{\mathrm{d}\epsilon}=\xi(\tilde{x},\tilde{t},\tilde{u},\tilde{v}),~~ \frac{\mathrm{d}\tilde{t}}{\mathrm{d}\epsilon}=\tau(\tilde{x},\tilde{t},\tilde{u},\tilde{v}),
~~\frac{\mathrm{d}\tilde{u}}{\mathrm{d}\epsilon}=\eta(\tilde{x},\tilde{t},\tilde{u},\tilde{v}),
~~\frac{\mathrm{d}\tilde{v}}{\mathrm{d}\epsilon}=\phi(\tilde{x},\tilde{t},\tilde{u},\tilde{v}),
\end{equation}
with the initial conditions
\begin{equation}\label{systemodes condition}
\tilde{x}(0)=x,~~\tilde{t}(0)=t,~~\tilde{u}(0)=u,~~\tilde{v}(0)=v.
\end{equation}

If $(u(x,t),v(x,t))$ is a solution of system \eqref{parabolic equations}, the actor of this system generated by $V$ can be denoted as
\begin{equation}\label{expuv}
(\tilde{u}_\epsilon(x,t), \tilde{v}_\epsilon(x,t))=\rho(\mathrm{exp}(\epsilon V))(u(x,t),v(x,t)),
\end{equation}
 where $\epsilon$ is the group parameter and $(\tilde{u}_\epsilon(x,t), \tilde{v}_\epsilon(x,t))$ is a new solution of system \eqref{parabolic equations}.
In addition,   in contrast to the scalar case, it is noticed that the fundamental solution is a $2\times2$ matrix for system \eqref{equation 2} from the definition of the fundamental solution of linear parabolic system in References \cite{KQ,KS}. Assume that $2\times2$ matrix $\bm{P}(t,x,y)$ is a solution of system \eqref{parabolic equations}, then
\begin{equation}\label{jibenjieshizi}
\bm{U}(x,t)=\int_R \bm{P}(t,x,y)\bm{f}(y)\mathrm{d}y,
\end{equation}
is a solution of the Cauchy problem for system \eqref{parabolic equations} with initial data $\bm{U}(x,0)=\bm{f}(x)$, where $\bm{U}=(u,v)^\mathrm{T}$. Next, by associating equations \eqref{expuv} and \eqref{jibenjieshizi}, we select two sets of independent solutions of system \eqref{parabolic equations}, denoted as $(u_i(x),v_i(x)),i=1,2$, which respond to two sets of group invariant solutions to system \eqref{parabolic equations}, denoted as $(\tilde{u}_\epsilon^i(x,t),\tilde{v}_\epsilon^i(x,t)), i=1,2$.

Furthermore, in order to apply the Laplace transform, sometimes we have to choose the proper transformation $\epsilon\rightarrow\lambda$ so that $(\tilde{u}_\epsilon^i(x,t),\tilde{v}_\epsilon^i(x,t))$ becomes $(\tilde{u}_\lambda^i(x,t),\tilde{v}_\lambda^i(x,t))$,  which satisfy the following conditions
\begin{equation}\label{uvchuzhi}
\tilde{u}^i_\lambda(x,0)=\mathrm{e}^{-\lambda x}L_i^u(x),~~\tilde{v}^i_\lambda(x,0)=\mathrm{e}^{-\lambda x}L_i^v(x),~~i=1,2.
\end{equation}

Finally, the following Theorem is given to present the formula of the fundamental solution of system \eqref{parabolic equations}.

\textbf{Theorem 3.1} For system \eqref{parabolic equations} with two dependent variables $(u,v)$ defined on $R^+\times[0,T]$. Suppose that its group invariant solutions $(\tilde{u}^i_\lambda(x,t),\tilde{v}^i_\lambda(x,t)), i=1,2$ are Laplace transformations in $y$. Namely, there exists a matrix
\begin{equation}\label{matrixabcd}
\bm{P}(t,x,y)=\left(\begin{matrix}
 A(t,x,y) & B(t,x,y) \\
  C(t,x,y) & D(t,x,y)
\end{matrix}\right)
\end{equation}
satisfying $$ \int_0^\infty \bm{P}(t,x,y)\bm{L}(y)\mathrm{e}^{-\lambda y}\mathrm{d}y=\bm{U}(x,t)$$
with the matrix
$$\bm{L}(x)=\left(\begin{matrix}
 L_1^u(x) & L_2^u(x) \\
  L_1^v(x) & L_2^v(x)
\end{matrix}\right),~~\bm{U}(x,t)=\left(\begin{matrix}
 \tilde{u}^1_\lambda(x,t) & \tilde{u}^2_\lambda(x,t) \\
  \tilde{v}^1_\lambda(x,t) & \tilde{v}^2_\lambda(x,t)
\end{matrix}\right),$$
then $\bm{P}(t,x,y)$ is the fundamental solution of system \eqref{parabolic equations}.

Since the proof of Theorem 3.1 is similar to the proof of Theorem 2.1 in Reference \cite{KQ}, we omit the details.
\subsection{Fundamental solution for system of PDEs}
\qquad It can be found that the key step to find the fundamental solution of system \eqref{parabolic equations} is to compute the expressions of $A(t,x,y)$, $B(t,x,y)$, $C(t,x,y)$, $D(t,x,y)$ in \eqref{matrixabcd} introduced in Subsection 3.1. In terms of Theorem 3.1, the fundamental solution can be constructed by taking the inverse Laplace transform of certain group invariant solution.
In the following, we present an example to illustrate the details.

\textbf{Example 3.1} Consider the system of variable coefficients parabolic equations
\begin{equation}\label{system Example}
\left\{
\begin{aligned}
&\mathcal{T}_t^\alpha u=xu_{xx}+av_x,\\
&\mathcal{T}_t^\alpha v=xv_{xx}+bu_x, ~~x>0,
\end{aligned}
\right.
\end{equation}
where $a$ and $b$ are constants and $ab>0.$

At the first step, from equations \eqref{vector fields v}-\eqref{yantuotiaojian}, we can obtain
\begin{equation}\label{systeminvariant condition}
\left\{
\begin{aligned}
&[\eta^{\alpha,t}-\xi u_{xx}-x\eta^{xx}-a\phi^x]|_{\eqref{system Example}}=0,\\
&[\phi^{\alpha,t}-\xi v_{xx}-x\phi^{xx}-b\eta^x]|_{\eqref{system Example}}=0.
\end{aligned}
\right.
\end{equation}
Equate the coefficients of $ u_x, u_{xx}, \cdots$ in system \eqref{systeminvariant condition} to be zero, which leads to the following determining equations
\begin{equation}\label{systemdetermining equations}
\left\{
  \begin{aligned}
  &t^{1-\alpha}\eta_t-x\eta_{xx}-a\phi_{x}=0,~~ x\bigg(-\tau_t+\frac{1-\alpha}{t}\tau\bigg)-\xi+2x\xi_x=0,\\
  &a\bigg(\eta_u-\tau_t+\frac{1-\alpha}{t}\tau\bigg)-2x\eta_{xv}-a(\phi_v-\xi_x)=0,\\
  &b\eta_v-t^{1-\alpha}\xi_t-x\bigg(2\eta_{xu}-\xi_{xx}\bigg)-a\phi_u=0,\\
  &t^{1-\alpha}\phi_t-x\phi_{xx}-b\eta_x=0,~~
  b\bigg(\phi_v-\tau_t+\frac{1-\alpha}{t}\tau\bigg)-2x\phi_{xu}-b(\eta_u-\xi_x)=0,\\
  &a\phi_u-t^{1-\alpha}\xi_t-x(2\phi_{xv}-\xi_{xx})-b\eta_v=0,\\
  &\tau_u=\tau_v=\tau_x=\xi_u=\xi_v=\eta_{uu}=\eta_{uv}=\eta_{vv}=\phi_{uu}=\phi_{uv}=\phi_{vv}.
  \end{aligned}
  \right.
\end{equation}
Solve equations \eqref{systemdetermining equations} to obtain a basis for Lie algebra of system \eqref{system Example}
\[V_1=t\partial_t+\alpha x\partial_x,~~   V_2=t^{1-\alpha}\partial_t,\]
\[V_3=t^{1+\alpha}\partial_t+2\alpha xt^\alpha\partial_x-(\alpha^2xu+a\alpha t^\alpha v)\partial_u-(\alpha^2xv+b\alpha t^\alpha u)\partial_v,\]
\[V_4=u\partial_u+v\partial_v,~~V_5=av\partial_u+bu\partial_v,~~V_{\phi_3}=\phi_3(x,t)\partial_v,~~V_{\eta_3}=\eta_3(x,t)\partial_u,\]
where $\eta_3$ and $\phi_3$ satisfy $t^{1-\alpha}\eta_{3t}-x\eta_{3xx}-a\phi_{3x}=0$ and $t^{1-\alpha}\phi_{3t}-x\phi_{3xx}-b\eta_{3x}=0$, respectively.

Now we are interested in vector field $V_3$ and intend to compute the group action generated by vector field $V_3$. Solve system \eqref{systemodes} with the initial conditions \eqref{systemodes condition} to yield
\begin{equation}\label{systemexampleuv}
\left\{
  \begin{aligned}
  \tilde{u}_\epsilon(x,t)=&\frac{1}{2\sqrt{ab}}\mathrm{e}^{-\frac{\alpha^2\epsilon x}{1+\alpha\epsilon t^\alpha}}\bigg(\bigg(\frac{\sqrt{ab}}{(1+\alpha\epsilon t^\alpha)^{\sqrt{ab}}}+\frac{\sqrt{ab}}{(1+\alpha\epsilon t^\alpha)^{-\sqrt{ab}}}\bigg)u\\
  &+\bigg(\frac{a}{(1+\alpha\epsilon t^\alpha)^{\sqrt{ab}}}-\frac{a}{(1+\alpha\epsilon t^\alpha)^{-\sqrt{ab}}}\bigg)v\bigg),\\
  \tilde{v}_\epsilon(x,t)=&\frac{1}{2a}\mathrm{e}^{-\frac{\alpha^2\epsilon x}{1+\alpha\epsilon t^\alpha}}\bigg(\bigg(\frac{\sqrt{ab}}{(1+\alpha\epsilon t^\alpha)^{\sqrt{ab}}}-\frac{\sqrt{ab}}{(1+\alpha\epsilon t^\alpha)^{-\sqrt{ab}}}\bigg)u\\
  &+\bigg(\frac{a}{(1+\alpha\epsilon t^\alpha)^{\sqrt{ab}}}+\frac{a}{(1+\alpha\epsilon t^\alpha)^{-\sqrt{ab}}}\bigg)v\bigg),
  \end{aligned}
  \right.
\end{equation}
where $u=u\big(\frac{x}{(1+\alpha\epsilon t^\alpha)^2},\frac{t}{(1+\alpha\epsilon t^\alpha)^{\frac1\alpha}}\big),$ $v=v\big(\frac{x}{(1+\alpha\epsilon t^\alpha)^2},\frac{t}{(1+\alpha\epsilon t^\alpha)^{\frac1\alpha}}\big).$ If $(u(x,t),v(x,t))$ solves system \eqref{system Example}, then $(\tilde{u}_\epsilon(x,t),\tilde{v}_\epsilon(x,t))$ is a new solution of system \eqref{system Example}.

Next, according to equation \eqref{matrixabcd}, we prove the existence of $\bm{P}(t,x,y)$, in other words, we need to show the explicit expressions for $A(t,x,y)$,  $B(t,x,y)$, $C(t,x,y)$ and $D(t,x,y)$. Choose two sets of solutions of system \eqref{system Example} of the following form
\begin{equation}\label{systemexampleuvchu}
(u_1,v_1)=\bigg(1,\frac{\sqrt{ab}}{a}\bigg),~~(u_2,v_2)=x^{1+\sqrt{ab}}\bigg(1,-\frac{\sqrt{ab}}{a}\bigg),
\end{equation}
then substitute  \eqref{systemexampleuvchu} into \eqref{systemexampleuv} to obtain
\begin{equation}\label{systemexampleuvchu1}
  \begin{aligned}
  &(\tilde{u}^1_\epsilon(x,t),\tilde{v}^1_\epsilon(x,t))=\bigg(\mathrm{e}^{-\frac{\alpha^2\epsilon x}{1+\alpha\epsilon t^\alpha}}\frac{1}{(1+\alpha\epsilon t^\alpha)^{\sqrt{ab}}},\frac{\sqrt{ab}}{a}\mathrm{e}^{-\frac{\alpha^2\epsilon x}{1+\alpha\epsilon t^\alpha}}\frac{1}{(1+\alpha\epsilon t^\alpha)^{\sqrt{ab}}}\bigg),\\
  &(\tilde{u}^2_\epsilon(x,t),\tilde{v}^2_\epsilon(x,t))=\bigg(\mathrm{e}^{-\frac{\alpha^2\epsilon x}{1+\alpha\epsilon t^\alpha}}\frac{x^{1+\sqrt{ab}}}{(1+\alpha\epsilon t^\alpha)^{\sqrt{ab}+2}},-\frac{\sqrt{ab}}{a}\mathrm{e}^{-\frac{\alpha^2\epsilon x}{1+\alpha\epsilon t^\alpha}}\frac{x^{1+\sqrt{ab}}}{(1+\alpha\epsilon t^\alpha)^{\sqrt{ab}+2}}\bigg),
  \end{aligned}
\end{equation}
which satisfies \[(\tilde{u}^1_\epsilon(x,0),\tilde{v}^1_\epsilon(x,0))=\mathrm{e}^{-\alpha^2\epsilon x}\bigg(1,\frac{\sqrt{ab}}{a}\bigg),~~(\tilde{u}^2_\epsilon(x,0),\tilde{v}^2_\epsilon(x,0))=\mathrm{e}^{-\alpha^2\epsilon x}x^{1+\sqrt{ab}}\bigg(1,-\frac{\sqrt{ab}}{a}\bigg).\]
In view of Theorem 3.1 and set $\lambda=\alpha^2\epsilon$ in  \eqref{systemexampleuvchu1}, we deduce that
\begin{equation}\label{AbCDshizi}
\left\{
  \begin{aligned}
  &\int_0^\infty(AL_1^u(y)+BL_1^v(y))\mathrm{e}^{-\lambda y}\mathrm{d}y=\mathrm{e}^{-\frac{\lambda x}{1+\lambda \frac{t^\alpha}{\alpha}}}\frac{1}{(1+\lambda \frac{t^\alpha}{\alpha})^{\sqrt{ab}}},\\
  &\int_0^\infty(AL_2^u(y)+BL_2^v(y))\mathrm{e}^{-\lambda y}\mathrm{d}y=\mathrm{e}^{-\frac{\lambda x}{1+\lambda \frac{t^\alpha}{\alpha}}}\frac{x^{1+\sqrt{ab}}}{(1+\lambda \frac{t^\alpha}{\alpha})^{\sqrt{ab}+2}},\\
    &\int_0^\infty(CL_1^u(y)+DL_1^v(y))\mathrm{e}^{-\lambda y}\mathrm{d}y=\frac{\sqrt{ab}}{a}\mathrm{e}^{-\frac{\lambda x}{1+\lambda \frac{t^\alpha}{\alpha}}}\frac{1}{(1+\lambda \frac{t^\alpha}{\alpha})^{\sqrt{ab}}},\\
  &\int_0^\infty(CL_2^u(y)+DL_2^v(y))\mathrm{e}^{-\lambda y}\mathrm{d}y=-\frac{\sqrt{ab}}{a}\mathrm{e}^{-\frac{\lambda x}{1+\lambda \frac{t^\alpha}{\alpha}}}\frac{x^{1+\sqrt{ab}}}{(1+\lambda \frac{t^\alpha}{\alpha})^{\sqrt{ab}+2}}.\\
  \end{aligned}
  \right.
\end{equation}
According to Lemma 2.2, we have
\begin{equation}\label{lapalceuv}
\left\{
  \begin{aligned}
  &\mathcal{L}\bigg(\mathrm{e}^{-\frac{\lambda x}{1+\lambda \frac{t^\alpha}{\alpha}}}\frac{1}{(1+\lambda \frac{t^\alpha}{\alpha})^{\sqrt{ab}}}\bigg)=\frac{\alpha}{t^\alpha}\mathrm{e}^{-\frac{\alpha(x+y)}{t^\alpha}}\bigg(\frac yx\bigg)^{\frac{\sqrt{ab}-1}{2}}I_{\sqrt{ab}-1}\bigg(\frac{2\alpha\sqrt{xy}}{t^\alpha}\bigg),\\
  &\mathcal{L}\bigg(\mathrm{e}^{-\frac{\lambda x}{1+\lambda \frac{t^\alpha}{\alpha}}}\frac{x^{1+\sqrt{ab}}}{(1+\lambda \frac{t^\alpha}{\alpha})^{\sqrt{ab}+2}}\bigg)=x^{1+\sqrt{ab}}\frac{\alpha}{t^\alpha}\mathrm{e}^{-\frac{\alpha(x+y)}{t^\alpha}}\bigg(\frac yx\bigg)^{\frac{1+\sqrt{ab}}{2}}I_{\sqrt{ab}+1}\bigg(\frac{2\alpha\sqrt{xy}}{t^\alpha}\bigg).\\
  \end{aligned}
  \right.
\end{equation}
Thanks to equations \eqref{lapalceuv} and inverting the Laplace transformation of equations \eqref{AbCDshizi} yields
\begin{align}\label{inverting Laplace ABCD1}
  &A+\frac{\sqrt{ab}}{a}B=\frac{\alpha}{t^\alpha}\mathrm{e}^{-\frac{\alpha(x+y)}{t^\alpha}}\bigg(\frac yx\bigg)^{\frac{\sqrt{ab}-1}{2}}I_{\sqrt{ab}-1}\bigg(\frac{2\alpha\sqrt{xy}}{t^\alpha}\bigg),
\end{align}
\begin{align}\label{inverting Laplace ABCD2}
  &A-\frac{\sqrt{ab}}{a}B=\frac{\alpha}{t^\alpha}\mathrm{e}^{-\frac{\alpha(x+y)}{t^\alpha}}\bigg(\frac yx\bigg)^{-\frac{1+\sqrt{ab}}{2}}I_{\sqrt{ab}+1}\bigg(\frac{2\alpha\sqrt{xy}}{t^\alpha}\bigg),
  \end{align}
\begin{align}\label{inverting Laplace ABCD3}
  &C+\frac{\sqrt{ab}}{a}D=\frac{\sqrt{ab}}{a}\frac{\alpha}{t^\alpha}\mathrm{e}^{-\frac{\alpha(x+y)}{t^\alpha}}\bigg(\frac yx\bigg)^{\frac{\sqrt{ab}-1}{2}}I_{\sqrt{ab}-1}\bigg(\frac{2\alpha\sqrt{xy}}{t^\alpha}\bigg),
  \end{align}
\begin{align}\label{inverting Laplace ABCD4}
  &C-\frac{\sqrt{ab}}{a}D=-\frac{\sqrt{ab}}{a}\frac{\alpha}{t^\alpha}\mathrm{e}^{-\frac{\alpha(x+y)}{t^\alpha}}\bigg(\frac yx\bigg)^{-\frac{1+\sqrt{ab}}{2}}I_{\sqrt{ab}+1}\bigg(\frac{2\alpha\sqrt{xy}}{t^\alpha}\bigg).
  \end{align}
Solving equations \eqref{inverting Laplace ABCD1}-\eqref{inverting Laplace ABCD4} and from \eqref{matrixabcd}, we obtain the following fundamental solution of system \eqref{system Example}
\begin{equation}\label{examplefundamental}
 \bm{P}(t,x,y)=\frac{\alpha}{2t^\alpha}\mathrm{e}^{-\frac{\alpha(x+y)}{t^\alpha}}\left(\begin{matrix}
 \gamma_1 & \frac{a}{\sqrt{ab}}\gamma_2 \\
  \frac{\sqrt{ab}}{a}\gamma_2 & \gamma_1
\end{matrix}\right),
\end{equation}
where \[\gamma_1=\bigg(\frac yx\bigg)^{\frac{\sqrt{ab}-1}{2}}I_{\sqrt{ab}-1}\bigg(\frac{2\alpha\sqrt{xy}}{t^\alpha}\bigg)+\bigg(\frac yx\bigg)^{-\frac{1+\sqrt{ab}}{2}}I_{\sqrt{ab}+1}\bigg(\frac{2\alpha\sqrt{xy}}{t^\alpha}\bigg),\]
\[\gamma_2=\bigg(\frac yx\bigg)^{\frac{\sqrt{ab}-1}{2}}I_{\sqrt{ab}-1}\bigg(\frac{2\alpha\sqrt{xy}}{t^\alpha}\bigg)-\bigg(\frac yx\bigg)^{-\frac{1+\sqrt{ab}}{2}}I_{\sqrt{ab}+1}\bigg(\frac{2\alpha\sqrt{xy}}{t^\alpha}\bigg).\]

This example shows that it is possible to derive the fundamental solution of conformable time fractional system \eqref{system Example} using the group invariant solution of the sysytem, so the question is whether we can obtain the fundamental solutions of other systems using the similar method as the one in Example 3.1.

Furthermore we discuss the following system
\begin{equation}\label{generalequation2}
\left\{
  \begin{aligned}
 &\mathcal{T}_t^\alpha u=x^mu_{xx}+a'v_x,\\
 &\mathcal{T}_t^\alpha v=x^mv_{xx}+b'u_x,
  \end{aligned}
  \right.
\end{equation}
which is more general than system \eqref{system Example}. If set $y=x^{\frac{2-q}{2}}$ and $\tau=(1-\frac{q}{2})^{\frac{2}{\alpha}} t$ in system \eqref{generalequation2}, this system  is transformed into the following conformable time fractional system
\begin{equation}\label{system Exampleset2}
\left\{
  \begin{aligned}
 &\mathcal{T}_\tau^\alpha u=u_{yy}+\frac{q}{(q-2)y}u_y+a'\frac{2}{2-q}y^{\frac{q}{q-2}}v_y,\\
 &\mathcal{T}_\tau^\alpha v=v_{yy}+\frac{q}{(q-2)y}v_y+b'\frac{2}{2-q}y^{\frac{q}{q-2}}u_y.
  \end{aligned}
  \right.
\end{equation}
In Section 3.3, we consider a more general conformable time fractional system \eqref{equation 2} than system \eqref{system Exampleset2}.
\subsection{Fundamental solution for system \eqref{equation 2}}
\qquad First, for the sake of simplicity, we consider only $mn>0$ in system \eqref{equation 2}. Recall that if the vector field $V=\xi(x,t,u,v)\partial_x+\tau(x,t,u,v)\partial_t+\eta(x,t,u,v)\partial_u+\phi(x,t,u,v)\partial_v$ generates a symmetry of system \eqref{equation 2}, then $V$ must satisfies
\begin{equation}\label{equation 2yantuotiaojian}
\left\{
\begin{aligned}
&\mathrm{Pr}^{(\alpha,2)}V\bigg(\mathcal{T}_t^\alpha u-u_{xx}-\frac{c}{x}u_x-mx^kv_x\bigg)|_{\eqref{equation 2}}=0,\\
&\mathrm{Pr}^{(\alpha,2)}V\bigg(\mathcal{T}_t^\alpha v-v_{xx}-\frac{c}{x}v_x-nx^ku_x\bigg)|_{\eqref{equation 2}}=0.
\end{aligned}
\right.
\end{equation}
Using the standard Lie point symmetry calculation algorithm and by means of equations \eqref{equation 2yantuotiaojian}, equating the coefficients of $u_x$, $u_{xx}, \cdots$ to be zero leads to
\begin{equation}\label{1equation 2}
t^{1-\alpha}\eta_t-\eta_{xx}-\frac{c}{x}\eta_x-mx^k\phi_x=0,
\end{equation}
\begin{equation}\label{2equation2}
t^{1-\alpha}\phi_t-\phi_{xx}-\frac{c}{x}\phi_x-nx^k\eta_x=0,
\end{equation}
\begin{equation}\label{3equation 2}
\frac{1-\alpha}{t}\tau-\tau_t+2\xi_x=0,
\end{equation}
\begin{equation}\label{4equation2}
\frac{c}{x}\bigg(-\tau_t+\frac{1-\alpha}{t}\tau\bigg)+nx^k\eta_v-t^{1-\alpha}\xi_t-(2\eta_{xu}-\xi_{xx})+\frac{c}{x^2}\xi+\frac{c}{x}\xi_x-mx^k\phi_u=0,
\end{equation}
\begin{equation}\label{5equation2}
\frac{c}{x}\bigg(-\tau_t+\frac{1-\alpha}{t}\tau\bigg)+mx^k\phi_u-t^{1-\alpha}\xi_t-(2\phi_{xv}-\xi_{xx})+\frac{c}{x^2}\xi+\frac{c}{x}\xi_x-nx^k\eta_v=0,
\end{equation}
\begin{equation}\label{6equation2}
mx^k\bigg(\eta_u-\tau_t+\frac{1-\alpha}{t}\tau\bigg)-2\eta_{xv}-mkx^{k-1}\xi-mx^k(\phi_v-\xi_x)=0,
\end{equation}
\begin{equation}\label{7equation2}
nx^k\bigg(\phi_v-\tau_t+\frac{1-\alpha}{t}\tau\bigg)-2\phi_{xu}-nkx^{k-1}\xi-nx^k(\eta_u-\xi_x)=0,
\end{equation}
\begin{equation}\label{8equation2}
\xi_u=\xi_v=\tau_x=\tau_u=\tau_v=\phi_{uu}=\phi_{uv}=\phi_{vv}=\eta_{uu}=\eta_{uv}=\eta_{vv}=0.
\end{equation}
Now it is time to solve equations \eqref{1equation 2}-\eqref{8equation2}. Consider equations \eqref{8equation2} to find
\begin{equation}\label{81equation2}
\tau=\tau(t),~~\xi=\xi_1(x,t),  \\
\end{equation}
and
\begin{equation}\label{811equation2}
\eta=\eta_1(x,t)u+\eta_2(x,t)v+\eta_3(x,t),~~ \phi=\phi_1(x,t)v+\phi_2(x,t)u+\phi_3(x,t). \\
\end{equation}
Substitute equations \eqref{81equation2}-\eqref{811equation2} into equations \eqref{1equation 2}-\eqref{7equation2} to obtain
\begin{equation}\label{11equation 2}
t^{1-\alpha}\eta_{1t}-\eta_{1xx}-\frac{c}{x}\eta_{1x}-mx^k\phi_{2x}=0,
\end{equation}
\begin{equation}\label{12equation 2}
t^{1-\alpha}\eta_{2t}-\eta_{2xx}-\frac{c}{x}\eta_{2x}-mx^k\phi_{1x}=0,
\end{equation}
\begin{equation}\label{13equation 2}
t^{1-\alpha}\eta_{3t}-\eta_{3xx}-\frac{c}{x}\eta_{3x}-mx^k\phi_{3x}=0,
\end{equation}
\begin{equation}\label{21equation2}
t^{1-\alpha}\phi_{1t}-\phi_{1xx}-\frac{c}{x}\phi_{1x}-nx^k\eta_{2x}=0,
\end{equation}
\begin{equation}\label{22equation2}
t^{1-\alpha}\phi_{2t}-\phi_{2xx}-\frac{c}{x}\phi_{2x}-nx^k\eta_{1x}=0,
\end{equation}
\begin{equation}\label{23equation2}
t^{1-\alpha}\phi_{3t}-\phi_{3xx}-\frac{c}{x}\phi_{3x}-nx^k\eta_{3x}=0,
\end{equation}
\begin{equation}\label{41equation2}
\frac{c}{x}\bigg(-\tau_t+\frac{1-\alpha}{t}\tau\bigg)+nx^k\eta_2-t^{1-\alpha}\xi_{1t}-2\eta_{1x}+\xi_{1xx}+\frac{c}{x^2}\xi_1+\frac{c}{x}\xi_{1x}-mx^k\phi_2=0,
\end{equation}
\begin{equation}\label{51equation2}
\frac{c}{x}\bigg(-\tau_t+\frac{1-\alpha}{t}\tau\bigg)+mx^k\phi_2-t^{1-\alpha}\xi_{1t}-2\phi_{1x}+\xi_{1xx}+\frac{c}{x^2}\xi_1+\frac{c}{x}\xi_{1x}-nx^k\eta_2=0,
\end{equation}
\begin{equation}\label{61equation2}
mx^k\bigg(\eta_1-\tau_t+\frac{1-\alpha}{t}\tau\bigg)-2\eta_{2x}-mkx^{k-1}\xi-mx^k(\phi_1-\xi_{1x})=0,
\end{equation}
\begin{equation}\label{71equation2}
nx^k\bigg(\phi_1-\tau_t+\frac{1-\alpha}{t}\tau\bigg)-2\phi_{2x}-nkx^{k-1}\xi-nx^k(\eta_1-\xi_{1x})=0,
\end{equation}
\begin{equation}\label{31equation 2}
\xi_1(x,t)=\frac{1}{2}\bigg(\tau_t-\frac{1-\alpha}{t}\tau\bigg)x+\sigma_1,
\end{equation}
where $\sigma_1=\sigma_1(t)$ is the undetermined function of $t$. Substituting \eqref{31equation 2} into  equations \eqref{41equation2}-\eqref{71equation2} and by addition or subtraction operation, we derive
\begin{equation}\label{4-5equation2}
x^k(n\eta_2-m\phi_2)+\phi_{1x}-\eta_{1x}=0,
\end{equation}
\begin{equation}\label{4+5equation2}
(-(1-\alpha)t^{-1-\alpha}\tau+(1-\alpha)t^{-\alpha}\tau_{t}-t^{1-\alpha}\tau_{tt})x-2t^{1-\alpha}\sigma_{1t}-2\eta_{1x}-2\phi_{1x}+\frac{2c}{x^2}\sigma_1=0,
\end{equation}
\begin{equation}\label{6-7equation2}
mnx^k(\eta_1-\phi_1)-n\eta_{2x}+m\phi_{2x}=0,
\end{equation}
\begin{equation}\label{6+7equation2}
2(m\phi_{2x}+n\eta_{2x})+mn(k+1)x^k\bigg(-\frac{1-\alpha}{t}\tau+\tau_t\bigg)+2mnkx^{k-1}\sigma_1=0,
\end{equation}
which lead to
\begin{equation}\label{phi2equation2}
\phi_2=\frac nm\eta_2+\frac{1}{mx^k}(\phi_{1x}-\eta_{1x}),
\end{equation}
\begin{equation}\label{phi1equation2}
\phi_1=\frac 14(-(1-\alpha)t^{-1-\alpha}\tau+(1-\alpha)t^{-\alpha}\tau_t-t^{1-\alpha}\tau_{tt})
x^2-t^{1-\alpha}\sigma_{1t}x-\eta_1-\frac{c}{x}\sigma_1+\sigma_2,
\end{equation}
where $\sigma_2=\sigma_2(t)$ is the undetermined function of $t$. Plug equation \eqref{phi1equation2} and equation \eqref{31equation 2} into equation \eqref{61equation2} and simplify them to obtain
\begin{align}\label{eta1equation2}
\eta_1=&\frac 18(-(1-\alpha)t^{-1-\alpha}\tau+(1-\alpha)t^{-\alpha}\tau_t-t^{1-\alpha}\tau_{tt})x^2-\frac 12t^{1-\alpha}\sigma_{1t}x\nonumber\\
&+\frac{k+1}{4}\bigg(\tau_t-\frac{1-\alpha}{t}\tau\bigg)+\frac{k-c}{2x}\sigma_1+\frac{1}{mx^k}\eta_{2x}+\frac12\sigma_2.
\end{align}
From equation \eqref{12equation 2}, we arrive at
\begin{align}\label{1-1equation 2}
&t^{1-\alpha}\eta_{2xt}-\frac{k+c}{x}\eta_{2xx}+\frac{k+c}{x^2}\eta_{2x}-
\frac14m(k+1)x^k(-(1-\alpha)t^{-1-\alpha}\tau\nonumber\\
&+(1-\alpha)t^{-\alpha}\tau_t-t^{1-\alpha}\tau_{tt})
+\frac12mkx^{k-1}t^{1-\alpha}\sigma_{1t}-\frac12m(k+c)(k-2)x^{k-3}\sigma_1=0.
\end{align}
Substituting equations \eqref{phi2equation2}-\eqref{eta1equation2} into equation \eqref{11equation 2} and equation \eqref{21equation2}  yields
\begin{align}\label{1-2equation 2}
&4x^{3-k}t^{1-\alpha}\eta_{2xt}-4(k+c)x^{2-k}\eta_{2xx}+4k(k+c)x^{1-k}\eta_{2x}
-m(k+1)x^3(-(1-\alpha)t^{-1-\alpha}\tau\nonumber\\
&+(1-\alpha)t^{-\alpha}\tau_t-t^{1-\alpha}\tau_{tt})+2mkx^2t^{1-\alpha}\sigma_{1t}
+2mk(k+c)\sigma_1=0.
\end{align}
Multiply equation \eqref{1-2equation 2} by $\frac14x^{k-3}$ and subtract equation \eqref{1-1equation 2} to get
\begin{align}\label{yizhongequation 2}
(k-1)(k+c)(m\sigma_1 x^{k-1}+\eta_{2x})=0.
\end{align}

Now we intend to provide Lie point symmetry admitted by system \eqref{equation 2}. To this end, we discuss it in two cases and in order to solve the fundamental solution and the conservation law of system \eqref{equation 2} later, we only consider the case $\eta_{2x}=0$.

\textbf{Case 3.1} $k+c\neq0$ and $k\neq1.$

There are two possibilities:

Subcase 3.1.1: $\tau\neq C_1t+C_2t^{1-\alpha}, C_1$ and $C_2$ are two arbitrary constants.

Thanks to equations \eqref{1-1equation 2}-\eqref{yizhongequation 2} and if $-(1-\alpha)t^{-1-\alpha}\tau+(1-\alpha)t^{-\alpha}\tau_t-t^{1-\alpha}\tau_{tt}\neq0$, we deduce that $k=-1$, $\sigma_1=0$, and $\eta_{2x}=0$. Furthermore, in this case we obtain the following vector fields
\[V_1=t\partial_t+\frac12\alpha x\partial_x, ~~V_2=t^{1-\alpha}\partial_t,\]
\begin{align*}
V_3=&t^{1+\alpha}\partial_t+\alpha xt^\alpha\partial_x-\bigg(\bigg(\frac{\alpha(c+1)t^\alpha}{2}+\frac{\alpha^2x^2}{4}\bigg)u+\frac{m\alpha t^\alpha}{2}v\bigg)\partial_u-\\
&\bigg(\bigg(\frac{\alpha(c+1)t^\alpha}{2}+\frac{\alpha^2x^2}{4}\bigg)v+\frac{n\alpha t^\alpha}{2}u\bigg)\partial_v,
\end{align*}
\[V_4=u\partial_u+v\partial_v,~~V_5=mv\partial_u+nu\partial_v,~~V_{\eta_3}=\eta_3\partial_u,~~V_{\phi_3}=\phi_3\partial_v.\]

Subcase 3.1.2: $\tau=C_1t+C_2t^{1-\alpha}.$

The basis for the Lie algebra is $V_{\eta_3}$, $V_{\phi_3}$, $V_2$, $V_4$ and $V_5$.

\textbf{Case 3.2} $k+c=0.$

In this case, we consider two subcases as follows:

Subcase 3.2.1: $k=-1$ $(c=1).$

The basis for the  Lie algebra of system \eqref{equation 2} is $V_{\eta_3}$, $V_{\phi_3}$, $V_1$, $V_2$, $V_3$, $V_4$, $V_5$.

Subcase 3.2.2: $k\neq-1.$

The basis for the  Lie algebra consists of  $V_{\eta_3}$, $V_{\phi_3}$, $V_2$,  $V_4$, $V_5$.

In the following, we use the Lie algebra obtained above to construct the fundamental solution of system \eqref{equation 2}.
Since the group action of  vector field $V_3$ is not trivial in $t$, let's consider vector field $V_3$,                                                                                                                                                                                                                                                                                                                                                                                                                                                                                                                                                                                                                                 which can be used to obtain the fundamental solution of system \eqref{equation 2} from the trivial solution of this system.

\textbf{Example 3.2} Consider the case $k=-1$ in system \eqref{equation 2}, namely
\begin{equation}\label{equation 3}
\left\{
\begin{aligned}
&\mathcal{T}_t^\alpha u=u_{xx}+\frac{c}{x}u_x+\frac{m}{x}v_x,\\
&\mathcal{T}_t^\alpha v=v_{xx}+\frac{c}{x}v_x+\frac{n}{x}u_x,~~ x>0.
\end{aligned}
\right.
\end{equation}
Due to the group action generated by $V_3$, we have the following result
\begin{equation}\label{eq2uvepsilon}
\left\{
  \begin{aligned}
  \tilde{u}_\epsilon(x,t)=&\frac{1}{2\sqrt{mn}}\mathrm{e}^{-\frac{\alpha^2\epsilon x^2}{4(1+\alpha\epsilon t^\alpha)}}\bigg(\bigg(\frac{\sqrt{mn}}{(1+\alpha\epsilon t^\alpha)^{\frac{c+1+\sqrt{mn}}{2}}}+\frac{\sqrt{mn}}{(1+\alpha\epsilon t^\alpha)^{\frac{c+1-\sqrt{mn}}{2}}}\bigg)u\\
  &+\bigg(\frac{m}{(1+\alpha\epsilon t^\alpha)^{\frac{c+1+\sqrt{mn}}{2}}}-\frac{m}{(1+\alpha\epsilon t^\alpha)^{\frac{c+1-\sqrt{mn}}{2}}}\bigg)v\bigg),\\
  \tilde{v}_\epsilon(x,t)=&\frac{1}{2m}\mathrm{e}^{-\frac{\alpha^2\epsilon x^2}{4(1+\alpha\epsilon t^\alpha)}}\bigg(\bigg(\frac{\sqrt{mn}}{(1+\alpha\epsilon t^\alpha)^{\frac{c+1+\sqrt{mn}}{2}}}-\frac{\sqrt{mn}}{(1+\alpha\epsilon t^\alpha)^{\frac{c+1-\sqrt{mn}}{2}}}\bigg)u\\
  &+\bigg(\frac{m}{(1+\alpha\epsilon t^\alpha)^{\frac{c+1+\sqrt{mn}}{2}}}+\frac{m}{(1+\alpha\epsilon t^\alpha)^{\frac{c+1-\sqrt{mn}}{2}}}\bigg)v\bigg),
  \end{aligned}
  \right.
\end{equation}
where $u=u\big(\frac{x}{1+\alpha\epsilon t^\alpha},\frac{t}{(1+\alpha\epsilon t^\alpha)^{\frac1\alpha}}\big),$ $v=v\big(\frac{x}{1+\alpha\epsilon t^\alpha},\frac{t}{(1+\alpha\epsilon t^\alpha)^{\frac1\alpha}}\big).$ If $(u,v)$ is a solution of system \eqref{equation 3}, then   equations \eqref{eq2uvepsilon} is also a solution of system \eqref{equation 3}. Here, we choose
\begin{equation}\label{system2euv12}
(u_1,v_1)=\bigg(1,\frac{\sqrt{mn}}{m}\bigg),~~(u_2,v_2)=x^{1+\sqrt{mn}-c}\bigg(-\frac{\sqrt{mn}}{n},1\bigg),
\end{equation}
which solve system \eqref{equation 3}.

Plug the above equations \eqref{system2euv12} into equations \eqref{eq2uvepsilon} and set $\lambda=\alpha^2\epsilon$ to obtain
\begin{equation}\label{system2uvchu}
  \begin{aligned}
  \bm{U}_\lambda(x,t)=\mathrm{e}^{-\frac{\lambda x^2}{(1+\frac{4t^\alpha}{\alpha}\lambda)}}\left(\begin{matrix}
 \frac{1}{(1+\frac{4t^\alpha}{\alpha}\lambda)^{\frac{c+1+\sqrt{mn}}{2}}} &  \frac{-mx^{1+\sqrt{mn}+c}}{\sqrt{mn}(1+\frac{4t^\alpha}{\alpha}\lambda)^{\frac{3+\sqrt{mn}-c}{2}}} \\
 \frac{\sqrt{mn}}{m(1+\frac{4t^\alpha}{\alpha}\lambda)^{\frac{c+1+\sqrt{mn}}{2}}} & \frac{ x^{1+\sqrt{mn}+c}}{(1+\frac{4t^\alpha}{\alpha}\lambda)^{\frac{3+\sqrt{mn}-c}{2}}}
\end{matrix}\right),
  \end{aligned}
\end{equation}
which satisfies \begin{equation}\label{system2uvchu1}
  \begin{aligned}
  \bm{U}_\lambda(x,0)=\mathrm{e}^{-\lambda x^2}\left(\begin{matrix}
 1 & -\frac{mx^{1+\sqrt{mn}-c}}{\sqrt{mn}} \\
  \frac{\sqrt{mn}}{m}& x^{1+\sqrt{mn}-c}
\end{matrix}\right).
  \end{aligned}
\end{equation}

In view of Theorem 3.1 and equation \eqref{system2uvchu1}, we arrive at
\begin{equation}\label{AbCDshizi11}
 \int_0^\infty(AL_1^u(y)+BL_1^v(y))\mathrm{e}^{-\lambda y^2}\mathrm{d}y=\frac{1}{(1+\frac{4t^\alpha}{\alpha}\lambda)^{\frac{c+1+\sqrt{mn}}{2}}}\mathrm{e}^{-\frac{\lambda x^2}{(1+\frac{4t^\alpha}{\alpha}\lambda)}},
 \end{equation}
\begin{equation}\label{AbCDshizi112}
  \int_0^\infty(AL_1^u(y)+BL_1^v(y))\mathrm{e}^{-\lambda y^2}\mathrm{d}y=\frac{-mx^{1+\sqrt{mn}+c}}{\sqrt{mn}(1+\frac{4t^\alpha}{\alpha}\lambda)^{\frac{3+\sqrt{mn}-c}{2}}}\mathrm{e}^{-\frac{\lambda x^2}{(1+\frac{4t^\alpha}{\alpha}\lambda)}},
  \end{equation}
  \begin{equation}\label{AbCDshizi113}
  \int_0^\infty(CL_2^u(y)+DL_2^v(y))\mathrm{e}^{-\lambda y^2}\mathrm{d}y=\frac{\sqrt{mn}}{m(1+\frac{4t^\alpha}{\alpha}\lambda)^{\frac{c+1+\sqrt{mn}}{2}}}\mathrm{e}^{-\frac{\lambda x^2}{(1+\frac{4t^\alpha}{\alpha}\lambda)}},
  \end{equation}
  \begin{equation}\label{AbCDshizi114}
  \int_0^\infty(CL_2^u(y)+DL_2^v(y))\mathrm{e}^{-\lambda y^2}\mathrm{d}y=\frac{x^{1+\sqrt{mn}+c}}{(1+\frac{4t^\alpha}{\alpha}\lambda)^{\frac{3+\sqrt{mn}-c}{2}}}\mathrm{e}^{-\frac{\lambda x^2}{(1+\frac{4t^\alpha}{\alpha}\lambda)}}.
\end{equation}
According to Lemma 2.2, we derive
\begin{equation}\label{lapalceuv1}
  \begin{aligned}
  &\mathcal{L}\bigg(\frac{1}{(1+\frac{4t^\alpha}{\alpha}\lambda)^{\frac{c+1+\sqrt{mn}}{2}}}\mathrm{e}^{-\frac{\lambda x^2}{(1+\frac{4t^\alpha}{\alpha}\lambda)}}\bigg)\\
  =&\frac{\alpha}{4t^\alpha}\mathrm{e}^{-\frac{\alpha(x^2+y^2)}{4t^\alpha}}\bigg(\frac yx\bigg)^{\frac{c+\sqrt{mn}-1}{2}}I_{\frac{c+\sqrt{mn}-1}{2}}\bigg(\frac{\alpha\sqrt{xy}}{2t^\alpha}\bigg),\\
  &\mathcal{L}\bigg(\frac{1}{(1+\frac{4t^\alpha}{\alpha}\lambda)^{\frac{3+\sqrt{mn}-c}{2}}}\mathrm{e}^{-\frac{\lambda x^2}{(1+\frac{4t^\alpha}{\alpha}\lambda)}}\bigg)\\
  =&x^{1+\sqrt{mn}-c}\frac{\alpha}{4t^\alpha}
  \mathrm{e}^{-\frac{\alpha(x^2+y^2)}{t^\alpha}}\bigg(\frac yx\bigg)^{\frac{1+\sqrt{mn}-c}{2}}I_{\frac{1+\sqrt{mn}-c}{2}}\bigg(\frac{\alpha\sqrt{xy}}{2t^\alpha}\bigg).\\
  \end{aligned}
\end{equation}
Inverting the Laplace transformation of equations \eqref{AbCDshizi11}-\eqref{AbCDshizi114} yields
\begin{equation}\label{inverting Laplace ABCD01}
  \frac{1}{2y}\bigg(A+\frac{\sqrt{mn}}{m}B\bigg)=\frac{\alpha}{4t^\alpha}\mathrm{e}^{-\frac{\alpha(x^2+y^2)}{4t^\alpha}}\bigg(\frac yx\bigg)^{\frac{c+\sqrt{mn}-1}{2}}I_{\frac{c+\sqrt{mn}-1}{2}}\bigg(\frac{\alpha\sqrt{xy}}{2t^\alpha}\bigg),
\end{equation}
\begin{align}\label{inverting Laplace ABCD02}
  &\frac{1}{2y}\bigg(-\frac{m}{\sqrt{mn}}A+B\bigg)y^{1+\sqrt{mn}-c}\nonumber\\
  ~~=&-\frac{mx^{1+\sqrt{mn}-c}}{\sqrt{mn}}\frac{\alpha}{4t^\alpha}
  \mathrm{e}^{-\frac{\alpha(x^2+y^2)}{t^\alpha}}\bigg(\frac yx\bigg)^{\frac{1+\sqrt{mn}-c}{2}}I_{\frac{1+\sqrt{mn}-c}{2}}\bigg(\frac{\alpha\sqrt{xy}}{2t^\alpha}\bigg),
\end{align}
\begin{align}\label{inverting Laplace ABCD03}
  \frac{1}{2y}\bigg(C+\frac{\sqrt{mn}}{m}D\bigg)=\frac{\sqrt{mn}}{m}\frac{\alpha}{4t^\alpha}\mathrm{e}^{-\frac{\alpha(x^2+y^2)}{4t^\alpha}}\bigg(\frac yx\bigg)^{\frac{c+\sqrt{mn}-1}{2}}I_{\frac{c+\sqrt{mn}-1}{2}}\bigg(\frac{\alpha\sqrt{xy}}{2t^\alpha}\bigg),
\end{align}
\begin{align}\label{inverting Laplace ABCD04}
  &\frac{1}{2y}\bigg(-\frac{m}{\sqrt{mn}}C+D\bigg)y^{1+\sqrt{mn}-c}\nonumber\\
  ~~=&\frac{\alpha x^{1+\sqrt{mn}-c}}{4t^\alpha}
  \mathrm{e}^{-\frac{\alpha(x^2+y^2)}{t^\alpha}}\bigg(\frac yx\bigg)^{\frac{1+\sqrt{mn}-c}{2}}I_{\frac{1+\sqrt{mn}-c}{2}}\bigg(\frac{\alpha\sqrt{xy}}{2t^\alpha}\bigg).
\end{align}

In view of equation \eqref{matrixabcd} and solve equations \eqref{inverting Laplace ABCD01}-\eqref{inverting Laplace ABCD04} for $A$, $B$, $C$ and $D$  to obtain the fundamental solution of system \eqref{equation 3}
\begin{equation}\label{system2fundamental}
 \bm{P}(t,x,y)=\frac{\alpha}{4t^\alpha}\mathrm{e}^{-\frac{\alpha(x^2+y^2)}{4t^\alpha}}
 \sqrt{xy}\left(\begin{matrix}
 \gamma_1 & \frac{m}{\sqrt{mn}}\gamma_2 \\
  \frac{\sqrt{mn}}{m}\gamma_2 & \gamma_1
\end{matrix}\right),
\end{equation}
where \[\gamma_1=\bigg(\frac yx\bigg)^{\frac{c+\sqrt{mn}}{2}}I_{\frac{c+\sqrt{mn}-1}{2}}\bigg(\frac{\alpha xy}{2t^\alpha}\bigg)+\bigg(\frac yx\bigg)^{\frac{c-\sqrt{mn}}{2}}I_{\frac{1+\sqrt{mn}-c}{2}}\bigg(\frac{\alpha xy}{2t^\alpha}\bigg),\]
\[\gamma_2=\bigg(\frac yx\bigg)^{\frac{c+\sqrt{mn}}{2}}I_{\frac{c+\sqrt{mn}-1}{2}}\bigg(\frac{\alpha xy}{2t^\alpha}\bigg)-\bigg(\frac yx\bigg)^{\frac{c-\sqrt{mn}}{2}}I_{\frac{1+\sqrt{mn}-c}{2}}\bigg(\frac{\alpha xy}{2t^\alpha}\bigg).\]
\subsection{Equivalence transformations and fundamental solutions}
\qquad Now,  we explore the relationship between the fundamental solutions of two systems related to the equivalent transformation
\begin{equation}\label{equivalent transformation}
 \tilde{t}=Y(t),~~ \tilde{x}=X(x,t),~~ \bm{\tilde{U}}=\bm{F}(x,t)\bm{U}(x,t),
\end{equation}
where $\bm{F}(x,t)=\left(\begin{matrix}
 r_1(x,t) & r_2(x,t) \\
  s_2(x,t) & s_1(x,t)
\end{matrix}\right),$ $\bm{\tilde{U}}=(\tilde{u},\tilde{v})^\mathrm{T}.$ Clearly, the invertibility of transformation \eqref{equivalent transformation} implies $X_x\neq0,$ $Y_t\neq0$ and $r_1s_1-r_2s_2\neq0.$

Write the inverse transformations of $X$ and $Y$ as
 $$x=Z(\tilde{x},\tilde{t}),~~t=Y^{-1}(\tilde{t}).$$
Assume that transformation \eqref{equivalent transformation} is an equivalence transformation of the class of conformable time fractional linear system
 \begin{equation}\label{initial equation}
 E(x,t,u,v)=0,~x\in\Omega,~t>0,
\end{equation}
so that the transformed system is
\begin{equation}\label{transformed equation}
 E(\tilde{x},\tilde{t},\tilde{u},\tilde{v})=0,~\tilde{x}\in\tilde{\Omega},~\tilde{t}>0,
\end{equation}
which belongs to the same class of system as the initial one.

If $\bm{U}(x,t)$ is a solution of initial system \eqref{initial equation} and from transformation \eqref{equivalent transformation}, then
\[\bm{\tilde{U}}(\tilde{x},\tilde{t})=\bm{F}(Z(\tilde{x},\tilde{t}),Y^{-1}(\tilde{t}))\bm{U}(Z(\tilde{x},\tilde{t}),
Y^{-1}(\tilde{t}))\]
is a solution of the transformed system \eqref{transformed equation}. Set $\tilde{t}(0)=0$ without loss of generality. In the following theorem, we show that if one has the fundamental solution to system \eqref{initial equation}, by the transformation \eqref{equivalent transformation}, one can get the fundamental solution to system \eqref{transformed equation}.

\textbf{Theorem 3.2} Assume that the linear system of PDEs \eqref{initial equation} can be transformed into system \eqref{transformed equation} by transformation \eqref{equivalent transformation} and the compatibility condition  $\tilde{t}(0)=0$ holds. If $\bm{\Gamma}(t,x,z)$ is the fundamental solution of system \eqref{initial equation}, then
\[\tilde{\bm{\Gamma}}(\tilde{t},\tilde{x},\tilde{z})=\bm{F}(Z(\tilde{x},\tilde{t}),\bm{Y}^{-1}(\tilde{t}))
\bm{\Gamma}(Y^{-1}(\tilde{t}),Z(\tilde{x},\tilde{t}),Z(\tilde{z},\tilde{t}))\bm{F}^{-1}(Z(\tilde{z},0),\tilde{t}(0))Z_{z}(\tilde{z},\tilde{t})\]
is a fundamental solution to the transformed system \eqref{transformed equation}.

\textbf{Proof} The proof of Theorem 3.2 is similar to the proof of Theorem 4.1 in Reference \cite{KQ}, thus in this paper, we omit it.

In the following, we consider the equivalence transformation for a class of linear conformable time fractional system
\begin{equation}\label{TFrac initial equation}
\left\{
\begin{aligned}
\mathcal{T}_t^\alpha u=h(x,t)u_{xx}+f_1(x,t)u_x+g_1(x,t)v_x,\\
\mathcal{T}_t^\alpha v=h(x,t)v_{xx}+f_2(x,t)v_x+g_2(x,t)u_x,
\end{aligned}
\right.
\end{equation}
in which system \eqref{system Example} and system \eqref{equation 2} are two special cases. Consider the invertible transformation
\begin{equation}\label{TFbianhuan}
\tilde{x}=X(x,t,u,v),~~\tilde{t}=Y(x,t,u,v),~~\tilde{u}=R(x,t,u,v),~~\tilde{v}=S(x,t,u,v),
\end{equation}
which preserves system \eqref{TFrac initial equation}. Namely, $\tilde{u}(\tilde{x},\tilde{t})$ and $\tilde{v}(\tilde{x},\tilde{t})$ satisfy system of the following form
\begin{equation}\label{TFrac transformed equation}
\left\{
\begin{aligned}
\mathcal{T}_{\tilde{t}}^\alpha \tilde{u}=h'(\tilde{x},\tilde{t})\tilde{u}_{\tilde{x}\tilde{x}}+f'_1(\tilde{x},\tilde{t})
\tilde{u}_{\tilde{x}}+g'_1(\tilde{x},\tilde{t})\tilde{v}_{\tilde{x}},\\
\mathcal{T}_{\tilde{t}}^\alpha \tilde{v}=h'(\tilde{x},\tilde{t})\tilde{v}_{\tilde{x}\tilde{x}}+f'_2(\tilde{x},\tilde{t})
\tilde{v}_{\tilde{x}}+g'_2(\tilde{x},\tilde{t})\tilde{u}_{\tilde{x}},
\end{aligned}
\right.
\end{equation}
where $h$, $h',$ $f_1$, $f'_1$, $f_2$, $f'_2$, $g_1$, $g'_1,$ $g_2$ and $g'_2$ are smooth functions of their arguments.

In view of transformation \eqref{TFbianhuan}, we can get expressions for $\mathcal{T}_t^\alpha u$, $\mathcal{T}_t^\alpha v$, $u_x$, $v_{x}$, $u_{xx}$ and $v_{xx}$. Since the transformation \eqref{TFbianhuan} preserves system \eqref{TFrac initial equation}, in other words, set the coefficients of $\tilde{u}_{\tilde{x}}\tilde{v}_{\tilde{t}}$, $\tilde{u}_{\tilde{t}}\tilde{v}_{\tilde{x}}$, $\tilde{u}_{\tilde{x}}^2$, $\tilde{u}_{\tilde{x}}\tilde{v}_{\tilde{x}}$, $\tilde{v}_{\tilde{x}}^2$ to be  zero, which leads to $T_x=T_u=T_v=0$ and $X_u=X_v=0,$ $ X_x\neq0$, $R_{uu}=R_{uv}=R_{vv}=S_{uu}=S_{uv}=S_{vv}=0$. Therefore, we can derive
\begin{equation}\label{TFbianhuanuvwan}
\left\{
\begin{aligned}
&\tilde{t}=Y(t),\\
&\tilde{x}=X(x,t),\\
&\tilde{u}=r_1(x,t)u+r_2(x,t)v+r_3(x,t),\\
&\tilde{v}=s_1(x,t)v+s_2(x,t)u+s_3(x,t),
\end{aligned}
\right.
\end{equation}
which implies
\begin{equation}\label{TFbianhuanuv}
\left\{
\begin{aligned}
&u=\frac{1}{\kappa}(s_1(x,t)\tilde{u}-r_2(x,t)\tilde{v}+\delta),\\
&v=\frac{1}{\kappa}(-s_2(x,t)\tilde{u}+r_1(x,t)\tilde{v}+\varrho),
\end{aligned}
\right.
\end{equation}
where $\kappa=r_1s_1-r_2s_2\neq0$, $\delta=r_2s_3-r_3s_1$, $\varrho=r_3s_2-r_1s_3$.
Consequently, according to transformation \eqref{TFbianhuanuvwan}, the expressions for $\mathcal{T}_t^\alpha u$, $\mathcal{T}_t^\alpha v$, $u_x$, $v_x$, $u_{xx}$ and $v_{xx}$ can be reduced to
\begin{equation}\label{TFbianhuanutvt1}
\begin{aligned}
\mathcal{T}_t^\alpha u=&\frac{1}{\kappa}[s_1\tilde{t}^{\alpha-1}\mathcal{T}_t^\alpha Y \mathcal{T}_{\tilde{t}}^\alpha \tilde{u}-r_2\tilde{t}^{\alpha-1}\mathcal{T}_t^\alpha Y \mathcal{T}_{\tilde{t}}^\alpha \tilde{v}+s_1\mathcal{T}_t^\alpha X\tilde{u}_{\tilde{x}}-r_2\mathcal{T}_t^\alpha X\tilde{v}_{\tilde{x}}\\
&+r_2(v\mathcal{T}_t^\alpha s_1+u\mathcal{T}_t^\alpha s_2+\mathcal{T}_t^\alpha s_3)-s_1(u\mathcal{T}_t^\alpha r_1+v\mathcal{T}_t^\alpha r_2+\mathcal{T}_t^\alpha r_3)],
\end{aligned}
\end{equation}
\begin{equation}\label{TFbianhuanutvt2}
\begin{aligned}
\mathcal{T}_t^\alpha v=&-\frac{1}{\kappa}[s_2\tilde{t}^{\alpha-1}\mathcal{T}_t^\alpha Y \mathcal{T}_{\tilde{t}}^\alpha \tilde{u}-r_1\tilde{t}^{\alpha-1}\mathcal{T}_t^\alpha Y \mathcal{T}_{\tilde{t}}^\alpha \tilde{v}+s_2\mathcal{T}_t^\alpha X\tilde{u}_{\tilde{x}}-r_1\mathcal{T}_t^\alpha X\tilde{v}_{\tilde{x}}\\
&+r_1(v\mathcal{T}_t^\alpha s_1+u\mathcal{T}_t^\alpha s_2+\mathcal{T}_t^\alpha s_3)-s_2(u\mathcal{T}_t^\alpha r_1+v\mathcal{T}_t^\alpha r_2+\mathcal{T}_t^\alpha r_3)],
\end{aligned}
\end{equation}
\begin{equation}\label{TFbianhuanuxvx1}
u_x=\frac{1}{\kappa}[s_1X_x\tilde{u}_{\tilde{x}}-r_2X_x\tilde{v}_{\tilde{x}}+r_2(s_{1x}v+s_{2x}u+s_{3x})
-s_1(r_{1x}u+r_{2x}v+r_{3x})],
\end{equation}
\begin{equation}\label{TFbianhuanuxvx2}
v_x=-\frac{1}{\kappa}[s_2X_x\tilde{u}_{\tilde{x}}-r_1X_x\tilde{v}_{\tilde{x}}+r_1(s_{1x}v+s_{2x}u+s_{3x})
-s_2(r_{1x}u+r_{2x}v+r_{3x})],
\end{equation}
\begin{equation}\label{TFbianhuanuxxvxx1}
\begin{aligned}
u_{xx}=&\frac{1}{\kappa}\{X_x^2(s_1\tilde{u}_{\tilde{x}\tilde{x}}-r_2\tilde{v}_{\tilde{x}\tilde{x}})
+X_{xx}(s_1\tilde{u}_{\tilde{x}}-r_2\tilde{v}_{\tilde{x}})+(s_{1xx}v+s_{2xx}u\\
&+s_{3xx})r_2-(r_{1xx}u+r_{2xx}v+r_{3xx})s_1
+\frac{1}{\kappa}[(s_1X_x\tilde{u}_{\tilde{x}}-r_2X_x\tilde{v}_{\tilde{x}}\\
&+(s_{1x}v+s_{2x}u+s_{3x})r_2-s_1(r_{1x}u+r_{2x}v+r_{3x}))(2s_{2x}r_2-2r_{1x}s_1)]\\
&+\frac{1}{\kappa}[(s_2X_x\tilde{u}_{\tilde{x}}
-r_1X_x\tilde{v}_{\tilde{x}}+(s_{1x}v+s_{2x}u+s_{3x})r_1-s_2(r_{1x}u+r_{2x}v\\
&+r_{3x}))(2s_{1x}r_2-2r_{2x}s_1)]\},
\end{aligned}
\end{equation}
\begin{equation}\label{TFbianhuanuxxvxx2}
\begin{aligned}
v_{xx}=&\frac{1}{\kappa}\{X_x^2(s_2\tilde{u}_{\tilde{x}\tilde{x}}-r_1\tilde{v}_{\tilde{x}\tilde{x}})
+X_{xx}(s_2\tilde{u}_{\tilde{x}}-r_1\tilde{v}_{\tilde{x}})+(s_{1xx}v+s_{2xx}u\\
&+s_{3xx})r_1-(r_{1xx}u+r_{2xx}v+r_{3xx})s_2
+\frac{1}{\kappa}[(s_1X_x\tilde{u}_{\tilde{x}}-r_2X_x\tilde{v}_{\tilde{x}}\\
&+(s_{1x}v+s_{2x}u+s_{3x})r_2-s_1(r_{1x}u+r_{2x}v+r_{3x}))(2s_{2x}r_1-2r_{1x}s_2)]\\
&+\frac{1}{\kappa}[(s_2X_x\tilde{u}_{\tilde{x}}
-r_1X_x\tilde{v}_{\tilde{x}}+(s_{1x}v+s_{2x}u+s_{3x})r_1-s_2(r_{1x}u+r_{2x}v\\
&+r_{3x}))(2s_{1x}r_1-2r_{2x}s_2)]\}.
\end{aligned}
\end{equation}
Next, substituting equations \eqref{TFbianhuanutvt1}-\eqref{TFbianhuanuxxvxx2} into system \eqref{TFrac initial equation} yields
\begin{align}\label{changshizi1}
&s_1\tilde{t}^{\alpha-1}\mathcal{T}_t^\alpha Y\mathcal{T}_{\tilde{t}}^\alpha \tilde{u}-r_2\tilde{t}^{\alpha-1}\mathcal{T}_t^\alpha Y\mathcal{T}_{\tilde{t}}^\alpha \tilde{v}-hs_1X_x^2\tilde{u}_{\tilde{x}\tilde{x}}+hr_2X_x^2\tilde{v}_{\tilde{x}\tilde{x}}
+\tilde{u}_{\tilde{x}}\big[s_1\mathcal{T}_t^\alpha X\nonumber\\&
-hs_1X_{xx}-2hX_xs_{1x}+\frac{2hs_1}{\kappa}X_x\kappa_x-f_1s_1X_x+g_1s_2X_x\big]+
\tilde{v}_{\tilde{x}}\big[-r_2\mathcal{T}_t^\alpha X\nonumber\\
&+hr_2X_{xx}+2hX_xr_{2x}-\frac{2hr_2}{\kappa}X_x\kappa_x+f_1r_2X_x-g_1r_1X_x\big]+
\tilde{u}\big[\kappa\mathcal{T}_t^\alpha (\frac{s_1}{\kappa})-h\kappa\big(\frac{s_1}{\kappa}\big)_{xx}\nonumber\\
&-f_1\kappa\big(\frac{s_1}{\kappa}\big)_x
+g_1\kappa\big(\frac{s_2}{\kappa}\big)_x\big]+\tilde{v}\big[-\kappa\mathcal{T}_t^\alpha (\frac{r_2}{\kappa})+h\kappa\big(\frac{r_2}{\kappa}\big)_{xx}+f_1\kappa\big(\frac{r_2}{\kappa}\big)_x
-g_1\kappa\big(\frac{r_1}{\kappa}\big)_x\big]\nonumber\\
&+\kappa\mathcal{T}_t^\alpha (\frac{\delta}{\kappa})-h\kappa\big(\frac{\delta}{\kappa}\big)_{xx}-f_1\kappa\big(\frac{\delta}{\kappa}\big)_x
-g_1\kappa\big(\frac{\varrho}{\kappa}\big)_x=0,
\end{align}
\begin{align}\label{changshizi2}
&-s_2\tilde{t}^{\alpha-1}\mathcal{T}_t^\alpha Y\mathcal{T}_{\tilde{t}}^\alpha \tilde{u}+r_1\tilde{t}^{\alpha-1}\mathcal{T}_t^\alpha Y\mathcal{T}_{\tilde{t}}^\alpha \tilde{v}+hs_2X_x^2\tilde{u}_{\tilde{x}\tilde{x}}-hr_1X_x^2\tilde{v}_{\tilde{x}\tilde{x}}
+\tilde{u}_{\tilde{x}}\big[-s_2\mathcal{T}_t^\alpha X\nonumber\\
&+hs_2X_{xx}+2hX_xs_{2x}-\frac{2hs_2}{\kappa}X_x\kappa_x+f_2s_2X_x-g_2s_1X_x\big]+
\tilde{v}_{\tilde{x}}\big[r_1\mathcal{T}_t^\alpha X-hr_1X_{xx}\nonumber\\
&-2hX_xr_{1x}+\frac{2hr_1}{\kappa}X_x\kappa_x-f_2r_1X_x+g_2r_2X_x\big]+
\tilde{u}\big[-\kappa\mathcal{T}_t^\alpha (\frac{s_2}{\kappa})+h\kappa\big(\frac{s_2}{\kappa}\big)_{xx}\nonumber\\
&+f_2\kappa\big(\frac{s_2}{\kappa}\big)_x
-g_2\kappa\big(\frac{s_1}{\kappa}\big)_x\big]+\tilde{v}\big[\kappa\mathcal{T}_t^\alpha (\frac{r_1}{\kappa})-h\kappa\big(\frac{r_1}{\kappa}\big)_{xx}-f_2\kappa\big(\frac{r_1}{\kappa}\big)_x
+g_2\kappa\big(\frac{r_2}{\kappa}\big)_x\big]\nonumber\\
&+\kappa\mathcal{T}_t^\alpha (\frac{\varrho}{\kappa})-h\kappa\big(\frac{\varrho}{\kappa}\big)_{xx}-f_2\kappa\big(\frac{\varrho}{\kappa}\big)_x
-g_2\kappa\big(\frac{\delta}{\kappa}\big)_x=0,
\end{align}
According to equations \eqref{changshizi1}-\eqref{changshizi2} and  transformed system \eqref{TFrac transformed equation}, we arrive at the following relations
\begin{align*}
h'(\tilde{x},\tilde{t})=\frac{h(Z(\tilde{x},\tilde{t}),Y^{-1}(\tilde{t}))X_x^2(Z(\tilde{x},\tilde{t}),Y^{-1}(\tilde{t}))}
{\tilde{t}^{\alpha-1}(\mathcal{T}_t^\alpha Y)(Y^{-1}(\tilde{t}))},
\end{align*}
\begin{align*}
f'_1(\tilde{x},\tilde{t})=&\frac{1}{\tilde{t}^{\alpha-1}(\mathcal{T}_t^\alpha Y)(Y^{-1}(\tilde{t}))}\bigg[-(\mathcal{T}_t^\alpha X)(Z(\tilde{x},\tilde{t}),Y^{-1}(\tilde{t}))+h(Z(\tilde{x},\tilde{t}),Y^{-1}(\tilde{t}))\\
&X_{xx}
(Z(\tilde{x},\tilde{t}),Y^{-1}(\tilde{t}))+\frac{2h(Z(\tilde{x},\tilde{t}),Y^{-1}(\tilde{t}))
X_x(Z(\tilde{x},\tilde{t}),Y^{-1}(\tilde{t}))}{\kappa}\big(r_{2x}s_2\\
&
-r_{1x}s_1\big)+
\frac{X_x(Z(\tilde{x},\tilde{t}),Y^{-1}(\tilde{t}))}{\kappa}
\big(f_1r_1s_1-g_1r_1s_2-f_2r_2s_2
+g_2r_2s_1\big)\bigg],
\end{align*}
\begin{align*}
g'_1(\tilde{x},\tilde{t})=&\frac{X_x(Z(\tilde{x},\tilde{t}),Y^{-1}(\tilde{t}))}{\kappa\tilde{t}^{\alpha-1}(\mathcal{T}_t^\alpha Y)(Y^{-1}(\tilde{t}))}\big[2h(Z(\tilde{x},\tilde{t}),Y^{-1}(\tilde{t}))(r_2
r_{1x}-r_1
r_{2x})-f_1r_1r_2\\
&+g_1r_1^2+f_2r_1
r_2-g_2r_2^2\big],
\end{align*}
\begin{align*}
f'_2(\tilde{x},\tilde{t})=&\frac{1}{\tilde{t}^{\alpha-1}(\mathcal{T}_t^\alpha Y)(Y^{-1}(\tilde{t}))}\bigg[-(\mathcal{T}_t^\alpha X)(Z(\tilde{x},\tilde{t}),Y^{-1}(\tilde{t}))+h(Z(\tilde{x},\tilde{t}),Y^{-1}(\tilde{t}))\\
&X_{xx}
(Z(\tilde{x},\tilde{t}),Y^{-1}(\tilde{t}))+\frac{2h(Z(\tilde{x},\tilde{t}),Y^{-1}(\tilde{t}))
X_x(Z(\tilde{x},\tilde{t}),Y^{-1}(\tilde{t}))}{\kappa}\big(r_{2}s_{2x}\\
&-r_{1}s_{1x}\big)+
\frac{X_x(Z(\tilde{x},\tilde{t}),Y^{-1}(\tilde{t}))}{\kappa}
\big(-f_1r_2s_2
+g_1r_1s_2
+f_2r_1s_1
-g_2r_2s_1\big)\bigg],
\end{align*}
\begin{align*}
g'_2(\tilde{x},\tilde{t})=&\frac{X_x(Z(\tilde{x},\tilde{t}),Y^{-1}(\tilde{t}))}{\kappa\tilde{t}^{\alpha-1}(\mathcal{T}_t^\alpha Y)(Y^{-1}(\tilde{t}))}\bigg[2h(Z(\tilde{x},\tilde{t}),Y^{-1}(\tilde{t}))(s_2
s_{1x}-s_1
s_{2x})+f_1s_1s_2
\\&-g_1s_2^2-f_2s_1
s_2+g_2
s_1^2\bigg],
\end{align*}
with $r_i,$ $s_i,$ $(i=1,2,3)$, satisfying
\[\mathcal{T}_t^\alpha(\frac{r_1}{\kappa})-h(\frac{r_1}{\kappa})_{xx}-f_2(\frac{r_1}{\kappa})_x
+g_2(\frac{r_2}{\kappa})_x=0,~~\mathcal{T}_t^\alpha(\frac{s_1}{\kappa})-h(\frac{s_1}{\kappa})_{xx}
-f_1(\frac{s_1}{\kappa})_x+g_1(\frac{s_2}{\kappa})_x=0,\]
\[\mathcal{T}_t^\alpha(\frac{r_2}{\kappa})-h(\frac{r_2}{\kappa})_{xx}-f_1(\frac{r_2}{\kappa})_x
+g_1(\frac{r_1}{\kappa})_x=0,~~\mathcal{T}_t^\alpha(\frac{s_2}{\kappa})-h(\frac{s_2}{\kappa})_{xx}
-f_2(\frac{s_2}{\kappa})_x+g_2(\frac{s_1}{\kappa})_x=0,\]
\[\mathcal{T}_t^\alpha(\frac{\delta}{\kappa})-h(\frac{\delta}{\kappa})_{xx}-f_1(\frac{\delta}{\kappa})_x
-g_1(\frac{\varrho}{\kappa})_x=0,~~\mathcal{T}_t^\alpha(\frac{\varrho}{\kappa})-h(\frac{\varrho}{\kappa})_{xx}
-f_2(\frac{\varrho}{\kappa})_x-g_2(\frac{\delta}{\kappa})_x=0.\]

\textbf{Example 3.3} System \eqref{equation 3} is related to the following system
\begin{equation}\label{TransformeExample}
\left\{
\begin{aligned}
&\mathcal{T}^\alpha_{\tilde{t}}\tilde{u}=\tilde{u}_{\tilde{x}\tilde{x}}-\bigg(c-2-\sqrt{mn}\frac{B_1}{B_2}\bigg)\frac{1}{\tilde{x}} \tilde{u}_{\tilde{x}}+\frac{2\sqrt{mn}a_1b_1}{B_2\tilde{x}}\tilde{v}_{\tilde{x}},\\
&\mathcal{T}^\alpha_{\tilde{t}}\tilde{v}=\tilde{v}_{\tilde{x}\tilde{x}}-\bigg(c-2+\sqrt{mn}\frac{B_1}{B_2}\bigg)\frac{1}{\tilde{x}} \tilde{v}_{\tilde{x}}-\frac{2\sqrt{mn}a_2b_2}{B_2\tilde{x}}\tilde{u}_{\tilde{x}},
\end{aligned}
\right.
\end{equation}
by the transformations
$$\left(\begin{matrix}
 \tilde{u}  \\
  \tilde{v}
\end{matrix}\right)=\bm{F}(x,t)\left(\begin{matrix}
 u  \\
  v
\end{matrix}\right),~~\tilde{x}=x,~~\tilde{t}=t,$$
where $$\bm{F}(x,t)=\left(\begin{matrix}
 \frac{1}{2B_2}(b_1x^{-\sqrt{mn}+c-1}-a_1x^{\sqrt{mn}+c-1}) & \frac{-\sqrt{mn}}{2nB_2}(a_1x^{\sqrt{mn}+c-1}+b_1x^{-\sqrt{mn}+c-1}) \\
  \frac{1}{2B_2}(a_2x^{\sqrt{mn}+c-1}-b_2x^{-\sqrt{mn}+c-1}) &\frac{\sqrt{mn}}{2nB_2}(a_2x^{\sqrt{mn}+c-1}+b_2x^{-\sqrt{mn}+c-1})
\end{matrix}\right),$$
$B_1=a_2b_1+a_1b_2, B_2=a_2b_1-a_1b_2\neq0.$ According to Theorem 3.2, the fundamental solutions of this system \eqref{TransformeExample} can be obtained.

\section{Conservation laws}
\qquad In this Section, we construct the conservation laws of the considered conformable fractional PDEs taking advantage of Lie algebras obtained above and  new Noether theorem \cite{I,IA}.

Consider the following conformable fractional differential equations
\begin{equation}\label{conservationsystem}
F_j(x,t,u_1,\cdots,u_s,\mathcal{T}_t^\alpha u_1,\cdots,\mathcal{T}_t^\alpha u_s,u_{1,x},\cdots, u_{s,x},\cdots)=0,~~j=1,\cdots,s,
\end{equation}
with two independent variables $(x,t)$ and $s(s>1)$ dependent variables $(u_1,\cdots,u_s)$. Assume that system \eqref{conservationsystem} admits the Lie symmetry generators written as follow
\begin{equation}\label{conservationsystemv}
V_i=\xi_i\partial_x+\tau_i\partial_t+\sum_{j=1}^s\eta^j_i\partial_{u_j},~~i=1,\cdots,n.
\end{equation}
The conserved vector $C=(C^t, C^x)$ for system \eqref{conservationsystem} satisfies the following conservation equation
\begin{equation}\label{L0}
\big(D_t(C^t)+D_x(C^x)\big)\mid_{\eqref{conservationsystem}}=0.
\end{equation}
The formal Lagrangian of system \eqref{conservationsystem} can be written as
\begin{equation}\label{LNS0}
L=\sum_{j=1}^{s}p_j(x,t)(F_j),
\end{equation}
with new dependent variable $ p_j(x,t)$, $j=1,\cdots,s$. The adjoint equations of formal Lagrangian \eqref{LNS0} are defined as follow \cite{IA}
\begin{equation}\label{CNSu}
F_j^*=\frac{\delta L}{\delta u_j}=0, ~~j=1,\cdots,s,
\end{equation}
where $\frac{\delta }{\delta u_j}$ is the Euler-Langrange operator denoted by
\begin{equation}\label{EL-CNS}
\frac{\delta}{\delta u_j}=\frac{\partial}{\partial u_j}+\sum_{l=1}^\infty (-1)^l D_{i1}D_{i2}\cdots D_{il}\frac{\partial}{\partial_{u_j,i_1i_2\cdots i_l}}.
\end{equation}
If the adjoint equations \eqref{CNSu} are satisfied for the solution of system \eqref{conservationsystem} upon the following substitutions
\begin{equation*}
p_j(x,t)=\psi_j(x,t,u_1,\cdots,u_s),
\end{equation*}
where $\psi_j\neq0$ for at least one $j$. It means that the following conditions must be held
\begin{equation*}
\frac{\delta L}{\delta u_j}|_{\eqref{conservationsystem}}=\sum_{i=1}^s\lambda_i^j(F_i).
\end{equation*}
For vector $V_i$, $i=1,\cdots,n$, conserved vectors can be obtained by the following formulas:
\begin{equation}\label{conservationsystemCi}
  \begin{aligned}
C^x_i=&\xi_i L+\sum_{j=1}^{s}\left(W^j_i\frac{\delta L}{\delta u_{j,x}}+D_x(W^j_i)\frac{\delta L}{\delta u_{j,xx}}+D^2_x(W^j_i)\frac{\delta L}{\delta u_{j,xxx}}+\cdots\right),\\
C^t_i=&\tau_i L+\sum_{j=1}^{s}\left(W^j_i\frac{\delta L}{\delta u_{j,t}}\right),
\end{aligned}
\end{equation}
where $W_i^j=\eta^j_i-\xi_iu_{j,x}-\tau_iu_{j,t},~i=1,\cdots,n, ~j=1,\cdots,s.$

\subsection{Conservation laws of system \eqref{system Example}}
\qquad Based on the symmetries admitted by system \eqref{system Example}, we intend to obtain the conservation law of system \eqref{system Example} in this Subsection.

The formal Lagrangian of system \eqref{system Example} is  written as
\begin{equation}\label{equation 1L00}
L=p(x,t)(\mathcal{T}_t^\alpha u-xu_{xx}-av_x)+
q(x,t)(\mathcal{T}_t^\alpha v-xv_{xx}-bu_x),
\end{equation}
with new dependent variable $ p(x,t) $ and $q(x,t)$.
The adjoint equations of formal Lagrangian equation \eqref{equation 1L00} are
\begin{equation}\label{equation 1F00}
\left\{\
\begin{aligned}
&\frac{\delta L}{\delta u}=F^*_1=-t^{1-\alpha}p_t-(1-\alpha)t^{-\alpha}p-2p_x+bq_x-xp_{xx}=0,\\
&\frac{\delta L}{\delta v}=F^*_2=-t^{1-\alpha}q_t-(1-\alpha)t^{-\alpha}q-2q_x+bp_x-xq_{xx}=0.
\end{aligned}
\right.
\end{equation}
Replace $ p(x,t)=\psi_1(x,t,u,v) $ and $ q(x,t)=\psi_2(x,t,u,v) $ in equations \eqref{equation 1F00} to derive
\begin{equation}\label{equation 1-p00}
\left\{\
\begin{aligned}
\frac{\delta L}{\delta u}|_{\{p=\psi_1\}}=\lambda_1(\mathcal{T}_t^\alpha u-xu_{xx}-av_x)+
\lambda_2(\mathcal{T}_t^\alpha v-xv_{xx}-bu_x),\\
\frac{\delta L}{\delta v}|_{\{q=\psi_2\}}=\lambda_3(\mathcal{T}_t^\alpha u-xu_{xx}-av_x)+
\lambda_4(\mathcal{T}_t^\alpha v-xv_{xx}-bu_x).
\end{aligned}
\right.
\end{equation}
According to equations \eqref{equation 1-p00}, we find that
\[\lambda_i=0 ~(i=1,2,3,4),~~\psi_1=(k_2+k_3x^{-1+\sqrt{ab}}+k_4x^{-1-\sqrt{ab}})t^{\alpha-1},\]
\[\psi_2=\frac{(k_1b+k_2+k_3x^{-1+\sqrt{ab}}\sqrt{ab}-k_4\sqrt{ab}x^{-1-\sqrt{ab}})t^{\alpha-1}}{b}.\]

Next, from the Lie algebra admitted by system \eqref{system Example}, by calculation, we obtain the following conserved vectors:\\
\textbf{Case 1} For $V_1=t\partial_t+\alpha x\partial_x,$ we can get
\begin{equation*}
 \begin{aligned}
C_1^x=&\frac{k_1}{b}\left(bxt^{\alpha} v_{xt}+ \alpha  b x^{2} t^{-1+\alpha}v_{xx}+\alpha b^{2} x t^{-1+\alpha}u_x +b^{2}u_tt^{\alpha} -b v_t t^{\alpha}\right)+\frac{k_2}{b}[x t^{\alpha} v_{xt}\\
&+\alpha  x^{2}t^{-1+\alpha} v_{xx}+ b x t^{\alpha} u_{xt}+ \alpha  b \,x^{2}t^{-1+\alpha} u_{xx}+\alpha abxt^{-1+\alpha} v_x+\alpha  bt^{-1+\alpha}u_x +(a b -1)t^{\alpha} v_t ]\\
&+\frac{k_3}{b}\big[\alpha t^{-1+\alpha}(\sqrt{a b} v_{xx}+b u_{xx}) x^{1+\sqrt{a b}}+((t^{\alpha} v_{xt}+\alpha  \,t^{-1+\alpha} v_x) \sqrt{a b}\\
&+b (t^{\alpha} u_{xt}+\alpha  \,t^{-1+\alpha} u_x)) x^{\sqrt{a b}}\big]+\frac{k_4}{b}[-\alpha t^{-1+\alpha}(\sqrt{a b}\, v_{xx}-b u_{xx}) x^{1-\sqrt{a b}}\\
&-((t^{\alpha} v_{xt}+\alpha  \,t^{-1+\alpha} v_x) \sqrt{a b}-b(t^{\alpha} u_{xt}+\alpha  \,t^{-1+\alpha} u_x)) x^{-\sqrt{a b}}],\\
C_1^t=&k_1\left(-\alpha  x v_x-t v_t\right)+\frac{k_2}{b}\left(-\alpha b  x u_x-b t u_t-\alpha  x v_x-t v_t\right) \\
&+\frac{k_3}{b}\left(-\alpha  \left(\sqrt{a b}\, v_x+bu_x  \right) x^{\sqrt{a b}}-\sqrt{a b}\, tx^{-1+\sqrt{a b}} v_t-b tx^{-1+\sqrt{a b}}u_t \right) \\
&+\frac{k_4}{b}\left(t \left(\sqrt{a b}\, v_t-bu_t \right) x^{-1-\sqrt{a b}}+\alpha  \left(\sqrt{a b}\, v_x-bu_x \right) x^{-\sqrt{a b}}\right).
 \end{aligned}
\end{equation*}
\textbf{Case 2} For $V_2=t^{1-\alpha}\partial_t,$ we have
\begin{equation*}
 \begin{aligned}
C_2^x=&k_1\left(x  v_{xt}+bu_t-v_t\right)+\frac{k_2}{b}\left(x v_{xt}+ b xu_{xt} +(a b -1) v_t\right)\\
&+\frac{k_3}{b}\left(\left(\sqrt{a b}\, v_{xt}+bu_{xt} \right) x^{\sqrt{a b}}\right)-\frac{k_4}{b} \left(\left(\sqrt{a b}\, v_{xt}-bu_{xt} \right) x^{-\sqrt{a b}}\right),\\
C_2^t=&-k_1t^{1-\alpha} v_t  +\frac{k_2}{b}\left((-bu_t -v_t) t^{1-\alpha}\right)-\frac{k_3}{b}\left(\left(\sqrt{a b}\, v_t+bu_t \right) x^{-1+\sqrt{a b}} t^{1-\alpha}\right)\\
&+\frac{k_4}{b} \left(\left(\sqrt{a b}\, v_t-bu_t \right) x^{-1-\sqrt{a b}} t^{1-\alpha}\right)
.
\end{aligned}
\end{equation*}
\textbf{Case 3} For $V_3=t^{1+\alpha}\partial_t+2\alpha xt^\alpha\partial_x-(\alpha^2xu+a\alpha t^\alpha v)\partial_u-(\alpha^2xv+b\alpha t^\alpha u)\partial_v,$ we derive that
 \begin{align*}
C_3^x=&\frac{k_1}{b}[(2\alpha b x^{2}v_{xx}+ \alpha b (3 b x u_x+abv -b u))t^{-1+2 \alpha}+bxt^{2 \alpha} v_{xt}+\alpha^2x (b x v_x+b^{2} u) t^{-1+\alpha}\\
&+(b^{2}u_t -b v_t) t^{2 \alpha}]+\frac{k_2}{b}[(2 x^{2} v_{xx}+ \alpha b (2 x^{2} u_{xx}+3 x u_x+3 a x v_x+(a b -1) u )) \,t^{-1+2 \alpha}\\
&+xt^{2 \alpha} v_{xt}+ b xt^{2 \alpha} u_{xt} + \alpha^{2}x (b x u_x+x v_x+b (a v+u)) t^{-1+\alpha}+(a b -1)t^{2 \alpha} v_t ]\\
&+\frac{k_3}{b}[ \alpha ((2 \sqrt{a b}\, v_{xx}+2 b u_{xx}) t^{-1+2 \alpha}+\alpha t^{-1+\alpha}   (\sqrt{a b}\, v_x+bu_x  ))x^{1+\sqrt{a b}}\\
&+(\alpha\sqrt{a b}(bu_x +2v_x )+\alpha b (a v_x+2 u_x)) x^{\sqrt{a b}}t^{-1+2 \alpha}+(\sqrt{a b}( \alpha^{2}t^{-1+\alpha}v +t^{2 \alpha} v_{xt})\\
&+b(\alpha^{2}t^{-1+\alpha}u+t^{2 \alpha} u_{xt})) x^{\sqrt{a b}}] +\frac{k_4}{b}[-\alpha  ((2 \sqrt{a b}\, v_{xx}-2 b u_{xx}) t^{-1+2 \alpha}\\
&+ \alpha t^{-1+\alpha}(\sqrt{a b}\, v_x-bu_x)) x^{1-\sqrt{a b}}+(\alpha  (\sqrt{a b}(-bu_x -2 v_x)+b (a v_x+2 u_x)) t^{-1+2 \alpha}\\
&+\sqrt{a b}(-\alpha^{2}t^{-1+\alpha}v-t^{2 \alpha} v_{xt})+b (\alpha^{2}t^{-1+\alpha}u +t^{2 \alpha} u_{xt})) x^{-\sqrt{a b}}] ,\\
C_3^t=&\frac{k_1}{b}[-bt^{1+\alpha}v_t-\alpha(2x bt^{\alpha}  v_x+b^{2} t^{\alpha}u +\alpha b x v)] +\frac{k_2}{b}[(-b u_t-v_t) t^{1+\alpha}-\alpha(2b x t^{\alpha}u_x\\
&+2 xt^{\alpha} v_x+b (a v+u ) t^{\alpha}+\alpha x (b u +v))]+\frac{k_3}{b}[-2 \alpha(\sqrt{a b}(t^{\alpha} v_x+\frac{\alpha v}{2}) \\
&+b (t^{\alpha} u_x+\frac{\alpha u}{2}))x^{\sqrt{a b}}-x^{-1+\sqrt{a b}} \sqrt{a b}(\alpha  bt^{\alpha} u +t^{1+\alpha} v_t)-bx^{-1+\sqrt{a b}}t^{1+\alpha} u_t \\
&-\alpha abx^{-1+\sqrt{a b}}t^{\alpha} v] +\frac{k_4}{b}[-(\sqrt{a b}(-\alpha bt^{\alpha} u -t^{1+\alpha} v_t)+b (\alpha at^{\alpha} v \\
&+t^{1+\alpha} u_t)) x^{-1-\sqrt{a b}}+2 \alpha(\sqrt{a b}(t^{\alpha} v_x+\frac{\alpha v}{2})-b (t^{\alpha} u_x+\frac{\alpha u}{2}))  \,x^{-\sqrt{a b}}].
\end{align*}
\textbf{Case 4} For $V_4=u\partial_u+v\partial_v,$ we arrive at
\begin{equation*}
 \begin{aligned}
C_4^x=&\frac{k_1}{b}\left(-b xv_x-b^{2}u+b v \right)t^{-1+\alpha}  +\frac{k_2}{b} \left(-x v_x-bx u_x+\left(-a b +1\right) v \right)t^{-1+\alpha} \\
&-\frac{k_3}{b}\left(\sqrt{a b}\, v_x+bu_x\right) x^{\sqrt{a b}}t^{-1+\alpha}  +\frac{k_4}{b}\left(\sqrt{a b} x^{-\sqrt{a b}} v_x-bx^{-\sqrt{a b}}u_x \right)t^{-1+\alpha} ,\\
C_4^t=&k_1v +\frac{k_2}{b}\left(b u +v\right)+\frac{k_3}{b}\left(\sqrt{a b}v+b u\right) x^{-1+\sqrt{a b}}-\frac{k_4}{b} \left(\sqrt{a b}v-b u\right) x^{-1-\sqrt{a b}}.
 \end{aligned}
\end{equation*}
\textbf{Case 5} For $V_5=av\partial_u+bu\partial_v,$ we find out
\begin{equation*}
 \begin{aligned}
C_5^x=&k_1 \left(bx u_x+abv -b u \right)t^{-1+\alpha} +k_2 \left(x u_x+a xv_x+\left(a b -1\right) u \right)t^{-1+\alpha}\\
 &+k_3 \left(a v_x+\sqrt{a b}u_x \right) x^{\sqrt{a b}}t^{-1+\alpha} +k_4 \left(a v_x-\sqrt{a b}u_x\right) x^{-\sqrt{a b}}t^{-1+\alpha}
,\\
C_5^t=&k_1bu+k_2\left(a v +u \right)+k_3\left(a v+\sqrt{a b}u \right) x^{-1+\sqrt{a b}}+k_4 \left(a v -\sqrt{a b}u\right) x^{-1-\sqrt{a b}}.
 \end{aligned}
\end{equation*}
\subsection{Conservation laws of system \eqref{equation 2}}
\qquad In this Subsection, we will construct the conservation law of system \eqref{equation 2}. For convenience, here considering $k=-1$, namely, we consider the conservation law of system \eqref{equation 3}.

Similar to the construction of conservation laws for system \eqref{system Example} and  based on  the Lie algebras obtained in Subsection 3.3, we deduce that the following conserved vectors:\\
\textbf{Case 1} For $V_1=t\partial_t+\frac12\alpha x\partial_x, $ we can get
 \begin{align*}
C_1^x=&\frac{k_1}{2 n}[2 nxt^{\alpha} v_{xt}+\alpha n \,x^{2} t^{-1+\alpha}v_{xx}+\alpha(c n v_x+ n^{2}u_x) xt^{-1+\alpha}-2 ((-c +1) n v_t \\
&-n^{2}u_t ) t^{\alpha}]+\frac{k_2}{2 n}[-2 (c -1) xt^{\alpha} v_{xt}-\alpha(c -1) x^{2} t^{-1+\alpha} v_{xx}+2n xt^{\alpha} u_{xt}  \\
&+ \alpha  n \,x^{2}t^{-1+\alpha} u_{xx}-\alpha((-m n +c (c -1)) v_x-nu_x) x  t^{-1+\alpha}-2 (-m n  \\
&+(c -1)^{2}) t^{\alpha}v_t]+\frac{k_3}{2 n}[\alpha  \,t^{-1+\alpha} (\sqrt{m n}\, v_{xx}+nu_{xx}) x^{c +\sqrt{m n}+1}+(\sqrt{m n}(\alpha  \,t^{-1+\alpha} v_x \\
&+2 t^{\alpha} v_{xt})+n (\alpha  \,t^{-1+\alpha} u_x+2 t^{\alpha} u_{xt})) x^{c +\sqrt{m n}}]+\frac{k_4}{2 n}[-\alpha  \,t^{-1+\alpha} (\sqrt{m n}\, v_{xx} \\
&-nu_{xx}) x^{c -\sqrt{m n}+1}-(\sqrt{m n}(\alpha  \,t^{-1+\alpha} v_x +2 t^{\alpha} v_{xt})-n (\alpha  \,t^{-1+\alpha} u_x+2 t^{\alpha} u_{xt})) x^{c -\sqrt{m n}}],\\
C_1^t=&-\frac{k_1}{2}\left(\alpha  x v_x+2 t v_t\right) x +\frac{k_2}{2 n}\left(-2nxtu_t -\alpha n x^{2} u_x+\left(c -1\right) \left(\alpha  x v_x+2 t v_t\right) x \right) \\
&+\frac{k_3}{2 n}\left(-2 t \left(\sqrt{m n}\, v_t+nu_t\right) x^{c +\sqrt{m n}}-\sqrt{m n}\, \alpha x^{c +\sqrt{m n}+1} v_x -n \alpha  \,x^{c +\sqrt{m n}+1} u_x\right) \\
&+\frac{k_4}{2 n}\left(\alpha  \left(\sqrt{m n}\, v_x-nu_x \right) x^{c -\sqrt{m n}+1}+2 t \left(\sqrt{m n}\, v_t-nu_t \right) x^{c -\sqrt{m n}}\right).
\end{align*}
\textbf{Case 2} For $V_2=t^{1-\alpha}\partial_t,$ we have
\begin{equation*}
 \begin{aligned}
C_2^x=&\frac{k_1}{n}\left(n x v_{xt}+n\left(c -1\right) v_t+n^{2}u_t\right) +\frac{k_2}{n}\left(\left(1-c\right) x v_{xt}+n xu_{xt}+\left(m n -\left(c -1\right)^{2}\right) v_t\right) \\
&+\frac{k_3}{n}\left(\sqrt{m n}\, v_{xt}+nu_{xt}\right) x^{c +\sqrt{m n}}-\frac{k_4}{n} \left(\sqrt{m n}\, v_{xt}-nu_{xt}\right) x^{c -\sqrt{m n}},\\
C_2^t=&-k_1xt^{1-\alpha} v_t +\frac{k_2}{n} \left(\left(c -1\right) v_t-nu_t \right) x t^{1-\alpha} -\frac{k_3}{n}\left(\sqrt{m n}\, v_t+nu_t\right) x^{c +\sqrt{m n}}t^{1-\alpha}  \\ &+\frac{k_4}{n}\left(\sqrt{m n}\, v_t-nu_t\right) x^{c -\sqrt{m n}}t^{1-\alpha}.
\end{aligned}
\end{equation*}
\textbf{Case 3} For \begin{align*}
V_3=&t^{1+\alpha}\partial_t+\alpha xt^\alpha\partial_x-\bigg(\bigg(\frac{\alpha(c+1)t^\alpha}{2}+\frac{\alpha^2x^2}{4}\bigg)u+\frac{m\alpha t^\alpha}{2}v\bigg)\partial_u-\\
&\bigg(\bigg(\frac{\alpha(c+1)t^\alpha}{2}+\frac{\alpha^2x^2}{4}\bigg)v+\frac{n\alpha t^\alpha}{2}u\bigg)\partial_v,
\end{align*} we obtain that
 \begin{align*}
C_3^x=&\frac{k_1}{4 n}[-2 \alpha(-2 n \,x^{2} v_{xx}-(1+3c) n x v_x-3 n^{2} x u_x+(-m \,n^{2}-(c^2 -1) n ) v \\
&-2 c \,n^{2}u)t^{-1+2 \alpha}+4 t^{2 \alpha} n x v_{xt}-\alpha^{2} (-n x v_x+(-c -1) n v- n^{2}u) x^{2} t^{-1+\alpha}\\
&-4 t^{2 \alpha} ((-c +1) n v_t-n^{2}u_t)]+\frac{k_2}{4 n}[-2 \alpha(2 (c -1) x^{2} v_{xx}-2 n \,x^{2} u_{xx}\\
&+(-3m n +(3c +1) (c -1)) x v_x-4 n x u_x+(-(c +1) m n +(c +1) (c -1)^{2}) v \\
&+u  (-m n +(c -1)^{2}) n )t^{-1+2 \alpha}-4 (c -1) x t^{2 \alpha} v_{xt}+4 n xt^{2 \alpha} u_{xt}  \\
&-\alpha^{2} ((c -1) x v_x-n x u_x+(c^{2}-m n -1) v-2 nu) x^{2} t^{-1+\alpha}-4 t^{2 \alpha} (-m n \\
&+(c -1)^{2}) v_t]+\frac{k_3}{4 n}[2\alpha ((2 \sqrt{m n}\, v_{xx}+2 nu_{xx}) t^{-1+2 \alpha}+\alpha  \,t^{-1+\alpha} (nu+\\
&\sqrt{m n}v))x^{c +\sqrt{m n}+1}+\alpha^{2} t^{-1+\alpha} (\sqrt{m n}\, v_x+nu_x) x^{c +\sqrt{m n}+2}+2 (\alpha ( \sqrt{m n}((c +3) v_x\\
&+nu_x)+n (m v_x+(c +3)u_x ))t^{-1+2 \alpha}+2 t^{2 \alpha} (\sqrt{m n}\, v_{xt}+nu_{xt})) x^{c +\sqrt{m n}}] \\
&+\frac{k_4}{4 n}[2\alpha((-2 \sqrt{m n}\, v_{xx}+2 nu_{xx}) t^{-1+2 \alpha}+\alpha  \,t^{-1+\alpha} (nu-\sqrt{m n}v))x^{c -\sqrt{m n}+1}\\
&-\alpha^{2} t^{-1+\alpha} (\sqrt{m n}\, v_x-nu_x) x^{c -\sqrt{m n}+2}-2 (\alpha(\sqrt{m n}((c +3) v_x+nu_x)-n (m v_x\\
&+(c +3)u_x)) t^{-1+2 \alpha}+2 t^{2 \alpha} (\sqrt{m n}\, v_{xt}-nu_{xt})) x^{c -\sqrt{m n}}],\\
C_3^t=&\frac{k_1}{2 n}[x (-2 n t^{1+\alpha}v_t +\alpha  (-2 x n \,t^{\alpha} v_x+((-c -1) n v- n^{2}u) t^{\alpha}-\frac{{\alpha}nx^{2} v} {2}))] \\
&+\frac{k_2}{2 n}[x ((2 (c -1) v_t-2 n u_t) t^{1+\alpha}+\alpha  (2 x (c -1) t^{\alpha} v_x-2 n xt^{\alpha} u_x +((c^{2}-m n -1) v\\
&-2 n u) t^{\alpha}+\frac{x^{2} ((c -1) v-nu) \alpha}{2}))]+\frac{k_3}{4 n}[-4 \alpha  \,t^{\alpha} (\sqrt{m n}\, v_x+nu_x) x^{c +\sqrt{m n}+1}\\
&-\alpha^{2} (nu+\sqrt{m n}v) x^{c +\sqrt{m n}+2}-2 (\sqrt{m n}(2 t^{1+\alpha} v_t+\alpha((c +1) v +nu)   t^{\alpha})\\
&+n (2 t^{1+\alpha} u_t+\alpha(m v +(c +1)u)t^{\alpha})) x^{c +\sqrt{m n}}]+\frac{k_4}{4 n}[4 \alpha  \,t^{\alpha} (\sqrt{m n}\, v_x-nu_x) x^{c -\sqrt{m n}+1}\\
&-\alpha^{2} (nu-\sqrt{m n}v) x^{c -\sqrt{m n}+2}-2 ( \sqrt{m n}(-2 t^{1+\alpha} v_t-\alpha((c +1) v +nu) t^{\alpha})\\
&+n (2 t^{1+\alpha} u_t+\alpha(m v +(c +1)u)t^{\alpha})) x^{c -\sqrt{m n}}].
\end{align*}
\textbf{Case 4} For $V_4=u\partial_u+v\partial_v,$ we arrive at
\begin{equation*}
 \begin{aligned}
C_4^x=&\frac{k_1}{n}(-x n v_x+(-c +1) n v-n^{2}u) t^{-1+\alpha} +\frac{k_2}{n}((c -1)x v_x- n xu_x +(-m n \\
&+(c -1)^{2}) v) t^{-1+\alpha}-\frac{k_3}{n}(\sqrt{m n}\, v_x+nu_x) x^{c +\sqrt{m n}} t^{-1+\alpha} \\
&+\frac{k_4}{n} (\sqrt{m n}\, v_x-nu_x) x^{c -\sqrt{m n}} t^{-1+\alpha},\\
C_4^t=&k_1xv -\frac{k_2}{n}\left(\left(c -1\right) v -nu \right)x +\frac{k_3}{n}\left(nu +\sqrt{m n}v\right) x^{c +\sqrt{m n}} +\frac{k_4}{n} \left(nu -\sqrt{m n}v\right) x^{c -\sqrt{m n}}.
 \end{aligned}
\end{equation*}
\textbf{Case 5} For $V_5=mv\partial_u+nu\partial_v,$ we find out
\begin{equation*}
 \begin{aligned}
C_5^x=&k_1(-nx u_x-(c -1) n u -m n v ) + k_2((c -1)x u_x-m xv_x  +(-m n +(c -1)^{2}) u )\\
 &-k_3 (\sqrt{m n}\, u_x+m v_x) x^{c +\sqrt{m n}}+k_4 (\sqrt{m n}\, u_x-m v_x) x^{c -\sqrt{m n}},\\
C_5^t=&k_1 n x u-k_2((c -1) u -m v )x +k_3 (\sqrt{m n}u+m v ) x^{c +\sqrt{m n}}\\
&-k_4 (\sqrt{m n}u -m v ) x^{c -\sqrt{m n}}.
 \end{aligned}
\end{equation*}
\section{Conclusions}
\qquad In this paper, we developed Lie symmetry method to construct the fundamental solution for the conformable time fractional system \eqref{equation 2} with variable coefficients. Firstly, in Example 3.1, we proved that it is possible to obtain the fundamental solutions to conformable time fractional system \eqref{system Example} associated with the group invariant solutions of this system and Laplace transform. Next, by considering a more general system \eqref{generalequation2} than system  \eqref{system Example}  and set transformation $y=x^{\frac{2-q}{2}}$, $\tau=(1-\frac{q}{2})^{\frac{2}{\alpha}} t$ in system \eqref{generalequation2} to yield system \eqref{system Exampleset2} and then a more general system \eqref{equation 2} than system \eqref{system Exampleset2} was considered. From the group action generated by the obtained nontrivial vector fields, we constructed the group invariant solutions of system \eqref{equation 2}. Then by two sets of steady-state solutions and inverting Laplace transform of group invariant solutions, the fundamental solutions of system \eqref{equation 2} with $k=-1$ were expressed in a matrix. And it is observed that the fundamental solution \eqref{examplefundamental} of system \eqref{system Example} and fundamental solution \eqref{system2fundamental} of system \eqref{equation 2} at $\alpha=1$ are exactly the same as the results obtained in Reference \cite{KQ}. In addition, we demonstrated that the fundamental solutions of two conformable fractional systems can be related by the equivalence transformation. Moreover, through Example 3.3, we can directly obtain the fundamental solution of system \eqref{TransformeExample} from the fundamental solution of system \eqref{equation 3} by equivalence transformation. Finally, the conservation laws of  systems  \eqref{system Example} and \eqref{equation 2} were derived by new Noether theorem.\\

\textbf{FUNDING}
This paper is submitted to the Special issue on Symmetry, Invariants, and their Applications in honor of Peter J. Olver's 70th Birthday. This work is supported by the National Natural Science Foundation of China (Grant No. 12271433)

\textbf{Conflict of Interest} The authors declare that they have no conflict of interest.

\pdfbookmark[1]{References}{ref}

\end{document}